\documentclass[12pt]{article}

\usepackage[dvipsnames]{xcolor}

\usepackage[margin={2cm,1.5cm}]{geometry}

\usepackage{graphicx}
\usepackage{float}
\usepackage{subcaption}

\usepackage{tikz}

\usepackage{listings}

\usepackage{multirow}

\usepackage[backend=bibtex,firstinits=true,sortcites,style=trad-plain]{biblatex}
\usepackage{enumerate}
\usepackage{amsmath,amsthm,amsfonts,bm}

\theoremstyle{plain}% default
\newtheorem{theorem}{Theorem}[section]
\newtheorem{proposition}[theorem]{Proposition}
\newtheorem{lemma}[theorem]{Lemma}
\newtheorem{corollary}[theorem]{Corollary}
\newtheorem{remark}[theorem]{Remark}

\theoremstyle{definition}

\newtheorem{hyp}{Hypothesis}

\usepackage{mathtools}
\newcommand{\R}{\mathbb{R}}

\newcommand{\buu}{\overline{\mathbf{u}}}

\newcommand{\bmu}{\overline{\mu}}

\newcommand{\bp}{\overline{p}}
\newcommand{\bphi}{\overline{\phi}}
\newcommand{\brho}{\overline{\rho}}

\newcommand{\nn}{\mathbf{n}}

\newcommand{\uu}{\mathbf{u}}

\newcommand{\JJ}{\mathbf{J}}
\newcommand{\DD}{\mathbf{D}}

\newcommand{\vK}{v_{K}}

\newcommand{\vL}{v_{L}}

\newcommand{\uKs}{u_{K^*}}

\newcommand{\Uh}{\mathcal{U}_h}
\newcommand{\Ph}{\mathcal{P}_h}

\usepackage[colorlinks=true]{hyperref}
\hypersetup{
linkcolor=BrickRed
,citecolor=Green
,filecolor=Mulberry
,urlcolor=NavyBlue
,menucolor=BrickRed
,runcolor=Mulberry
,linkbordercolor=BrickRed
,citebordercolor=Green
,filebordercolor=Mulberry
,urlbordercolor=NavyBlue
,menubordercolor=BrickRed
,runbordercolor=Mulberry
}

\def\escalar#1#2{\left(#1,#2\right)}

\def\escalarL#1#2{\escalar{#1}{#2}}
\def\escalarLd#1#2{\escalar{#1}{#2}}
\def\escalarML#1#2{\escalar{#1}{#2}_h}

\def\T{\mathcal{T}}
\def\E{\mathcal{E}}

\def\Ehi{\mathcal{E}_h^\text{i}}
\def\Ehb{\mathcal{E}_h^\text{b}}

\def\Pd{\mathbb{P}^{\text{disc}}}
\def\Pc{\mathbb{P}^{\text{cont}}}
\def\Pb{\mathbb{P}^{\text{bubble}}}
\def\Pih{\Pi^h}

\def\norma#1{\left\|#1\right\|}

\def\salto#1{\left[\!\left[#1\right]\!\right]}
\def\media#1{\left\{\!\!\left\{#1\right\}\!\!\right\}}

\def\th#1#2#3#4#5#6{s_h^1(#1,#2,#3,#4,#5,#6)}

\def\sh#1#2#3#4{s_h^2(#1,#2,#3,#4)}
\def\shd#1#2#3#4{s_h^{2,\delta}(#1,#2,#3,#4)}

\def\ch#1#2#3{c_h(#1,#2,#3)}

\def\aupw#1#2#3{a_h^{\text{upw}}(#1;#2,#3)}
\def\bupw#1#2#3{b_h^{\text{upw}}(#1;#2,#3)}

\def\sign#1{\text{sign}(#1)}

\title{\textbf{Property-preserving numerical approximation of a Cahn--Hilliard--Navier--Stokes model with variable density and degenerate mobility}}
\makeatletter
\def\@fnsymbol#1{\ensuremath{\ifcase#1\or *\or\dagger\or \ddagger\or \mathsection \or\mathparagraph\or *\or **\or \dagger\dagger \or \ddagger\ddagger \else\@ctrerr\fi}}
\makeatother

\author{Daniel Acosta-Soba\thanks{Departamento de Matemáticas, Universidad de Cádiz, Puerto Real, 11510 Cádiz, Spain -- Email: \texttt{daniel.acosta@uca.es} -- Corresponding author} \thanks{Department of Mathematics, University of Tennessee at Chattanooga, Chattanooga, TN 37403, USA}~,
~Francisco Guillén-González\thanks{Departamento de Ecuaciones Diferenciales y Análisis Numérico \& IMUS, Universidad de Sevilla, 41012 Seville, Spain -- Email: \texttt{guillen@us.es}}~,\\
~J. Rafael Rodríguez-Galván\thanks{Departamento de Matemáticas, Universidad de Cádiz, Puerto Real, 11510 Cádiz, Spain -- Email: \texttt{rafael.rodriguez@uca.es}}~,
~Jin Wang \thanks{Department of Mathematics, University of Tennessee at Chattanooga, Chattanooga, TN 37403, USA  -- Email: \texttt{jin-wang02@utc.edu}}}
\date{}

\begin{document}

\maketitle

\begin{abstract}
	 In this paper, we present a new computational framework to approximate a Cahn--Hilliard--Navier--Stokes model with variable density and degenerate mobility that preserves the mass of the mixture, the pointwise bounds of the density and the decreasing energy.
	This numerical scheme is based on a finite element approximation for the 
	 Navier--Stokes fluid flow with discontinuous pressure and an upwind discontinuous Galerkin scheme for the Cahn--Hilliard part.
	Finally, several numerical experiments such as a convergence test and some well-known benchmark problems are conducted.
\end{abstract}

\paragraph{Keywords:} Mass-conservation.  Discrete pointwise bounds. Discrete energy stability. 
Finite elements. Discontinuous Galerkin. Upwind scheme.

\section{Introduction}

Hydrodynamics has been considered a research field of increasing interest among the scientific community during the last few decades. In this sense, diffuse interface models were proposed as a successful alternative to model fluid-solid interaction after van der Waals introduced the foundations in the pioneering paper \cite{van1879thermodynamic}. Afterwards, these ideas were extended to fluid mixture and several works were published in this regard. In particular, both Hohelberg and Halpering, \cite{hohenberg1977theory}, and Gurtin et al., \cite{gurtin1996two}, arrived by different approaches to the same model, the well-known \textit{Model H}, which would lead to the 
Cahn--Hilliard--Navier--Stokes (CHNS) system.

Since then, many different CHNS models have been developed using different techniques and extended to the case of fluids with different densities, see the model by Boyer \cite{boyer2002theoretical} or by Ding et al. \cite{ding2007diffuse}. Moreover, several of these recent models satisfy 
some laws of thermodynamics.
This is the case for the model by Lowengrub and Truskinovsky, \cite{lowengrub1998quasi}, or the one by Abels et al., \cite{abels_thermodynamically_2011}, which introduces an extra convective term in the momentum equation due to the different densities of the fluids. In \cite{kim_2012} a careful revision of several CHNS models and their applications is provided. Also, recently, a very interesting survey has been published, \cite{ten2023unified}, in which the authors, Eikelder et al., discuss different existing well-known CHNS models analyzing their advantages and disadvantages from a physical point of view. In fact, the authors of \cite{ten2023unified} provide some notions on 
properties a CHNS model has to satisfy in order to be physically consistent.

One characteristic that many of these models share is that the density of the mixture is usually interpolated as a linear function of the phase-field function. Hence, ensuring the pointwise bounds for this phase-field function in the Cahn-Hilliard equation, for instance, by using a degenerate mobility (see \cite{acosta-soba_upwind_2022}) is crucial to ensure a physically consistent model. Also, CHNS models conserve the total mass of the mixture and, as mentioned above, they tend to be thermodynamically consistent in the sense that the solutions of these models usually minimize an underlying energy law. Therefore, as these properties are extremely important for the physical meaning of the models it is likewise important to preserve them when approximating their solutions.

However, the transport of the diffuse interface by the velocity of the fluid is typically modeled by means of a convective term that is introduced into the Cahn-Hilliard equation and, as shown in previous studies such as \cite{acosta-soba_upwind_2022}, this term may lead to numerical instabilities in highly convective regimes if it is not treated carefully. The instabilities result in nonphysical spurious oscillations that make the approximation of the phase-field variable lose the pointwise bounds. In this regard, removing the numerical instabilities in the case of the convective Cahn-Hilliard model has been an object of study in several recent works, see \cite{frank2018finite} or \cite{acosta-soba_upwind_2022}, where in the latter the authors enforce the pointwise bounds by means of a discontinuous Galerkin (DG) upwind technique. Different ideas such as the use of limiters have been used in the case of the CHNS systems. For instance, in \cite{liu2022pressure}, the authors 
developed, by means of flux and slope limiters, a bound-preserving decoupled approximation of a CHNS simplified system with constant mobility. Later, the same model was approximated  
 by high order polynomials using a decoupled scheme and a convex optimization technique with a scaling limiter to ensure the pointwise bounds, see \cite{liu2023simple}. In this line, the recent work \cite{guillentierra2024} has presented a numerical approximation of a CHNS which is mass conservative, energy stable and approximately pointwise bounded.

In addition, designing an approximation that satisfies a discrete version of the continuous energy in the diffuse-interface models is not straightforward and usually requires the use of specific time-discrete approximations such as the standard convex-splitting technique, \cite{eyre_1998_unconditionally}, or the more recently developed SAV approach, \cite{shen2018scalar}.
In this sense, several advancements have been made towards the approximation of the CHNS models preserving the energy-stability constraint. For instance, we can find the work \cite{tierra_guillen_abels_2014} where the authors propose an approximation of the model in \cite{abels_thermodynamically_2011} that decouples the phase-field equations from the fluid equations through a modified velocity. This approach was further studied in \cite{grun_guillen-gonzalez_metzger_2016} and extended to a fully decoupled approximation that uses a pressure correction approach, \cite{shen2015decoupled}. Other fractional time-stepping energy-stable discretizations of CHNS models can be found in \cite{salgado2013diffuse,deteix2022new,liu2015decoupled}.

Nevertheless, although it has been achieved in the case of a CHNS with a Flory-Huggins logarithmic potential (see \cite{chen2022positivity}), to our best knowledge there is no published work on an approximation of a CHNS model with a Ginzburg-Landau polynomial potential and degenerate mobility that ensures both the mass-conservation, pointwise bounds and energy-stability properties.

To address this challenge,
in this work, we provide an upwind DG approximation of the model by Abels et al. \cite{abels_thermodynamically_2011} where all the mass-conservation, the pointwise bounds and the energy-stability properties are preserved. 

Firstly,
in Section~\ref{sec:model} we introduce the CHNS model that we are going to consider and we present its properties.
Then, in Section~\ref{sec:coupled_scheme} we develop the structure-preserving approximation of the aforementioned model, showing that it satisfies all the mass-conservation, pointwise bounds and energy-stability properties. Finally, in Section~\ref{sec:numerical_experiments} we conduct several numerical experiments. First, we compute a preliminary accuracy test in Subsection~\ref{test:accuracy} for all the variables in both $L^2(\Omega)$ and $H^1(\Omega)$ norms. Then, we provide a simple test where two bubbles are mixed in Subsection~\ref{test:circle}. The results are in accordance with the previous theoretical analysis. Finally, in Subsections~\ref{test:bubble} and \ref{test:rayleigh} we couple the CHNS system with a term modeling the action of gravitational forces and conduct two benchmark tests: a heavier bubble in a lighter medium and a Rayleigh-Taylor type instability.

\section{Cahn--Hilliard--Navier--Stokes model}
\label{sec:model}

Let $\Omega\subset \R^d$ be a bounded polygonal domain. We consider a mixture of two fluids with different densities $0<\rho_1<\rho_2$ and introduce a phase-field function $\phi=\phi(x,t)\in [-1,1]$ such   that $\phi=-1 $ corresponds with fluid of density $\rho_1$, $\phi=1 $ with fluid of density $\rho_2$ and $\phi\in (-1,1)$ in the interface between the two fluids. Then, 
the diffuse-interface Cahn--Hilliard--Navier--Stokes model proposed by Abels et al. in \cite{abels_thermodynamically_2011} and further numerically studied in \cite{tierra_guillen_abels_2014,grun_guillen-gonzalez_metzger_2016,shen2015decoupled}, can be written as follows: 
\begin{subequations}
	\label{NS-CH_model}
	\begin{align}
		\label{eq:NS-CH_model_u}
		\rho(\phi)\uu_t+\left((\rho(\phi)\uu - \JJ
		)\cdot\nabla\right)\uu-\nabla\cdot(2\eta(\phi)\DD\uu)+\nabla p 
		+ \phi \nabla \mu&=0\quad\text{in }\Omega\times(0,T),\\
		\label{eq:NS-CH_model_p}
		\nabla\cdot\uu&=0\quad\text{in }\Omega\times(0,T),\\
		\label{eq:NS-CH_model_phi}
		\phi_t+\nabla\cdot(\phi\uu)-\nabla\cdot(M(\phi)\nabla \mu)&=0\quad\text{in }\Omega\times(0,T),\\
		\label{eq:NS-CH_model_mu}
		-\lambda\varepsilon\Delta\phi+\frac{\lambda}{\varepsilon}f(\phi)&=\mu\quad\text{in }\Omega\times(0,T),\\
		\label{eq:NS-CH_model_ic}
		\uu(0)=\uu_0,\quad\phi(0)&=\phi_0\quad\text{in }\Omega,\\
		\label{eq:NS-CH_model_bc}
		\uu=0,\quad\nabla\phi\cdot\nn=0,\quad M(\phi)\nabla\mu\cdot\nn&=0\quad\text{on }\partial\Omega.
	\end{align}
\end{subequations}
Here, $\uu$ and $p$ are the mean velocity and the pressure of the fluid respectively, and $\mu$ is the chemical potential related to the phase-field function $\phi$.
Also, $\DD\uu = \frac{1}{2}(\nabla\uu+\nabla\uu^t)$ is the strain tensor, $f(\phi)$ is the derivative of the Ginzburg-Landau double well potential $ F(\phi)=\frac{1}{4}(\phi^2-1)^2$, i.e. $f(\phi)=F'(\phi)=(\phi^2-1)\phi$, $M(\phi)=(1-\phi^2)_\oplus$ is the degenerate (truncated) mobility function and
$$\JJ=\frac{\rho_2-\rho_1}{2}M(\phi)\nabla\mu
$$ is the extra-convective term due to different densities. 
Moreover, the density of the mixture $\rho=\rho(\phi)$ depending on the phase-field variable $\phi$, can be defined   either as the solution of the mass balance equation
\begin{equation}
	\label{eq:NS-CH_model_rho}
	\escalarL{\partial_t\rho}{\overline\rho}-\escalarL{\rho\uu-\JJ}{\nabla\overline\rho}=0,\quad\forall\overline\rho\in H^1(\Omega), 
\quad  \hbox{in $(0,T)$},	
\end{equation}
 or, by taking into account the equation \eqref{eq:NS-CH_model_phi}, as the explicit relation
\begin{equation}
	\label{rho}
	\rho(\phi)=\frac{\rho_1+\rho_2}{2}+\frac{\rho_2-\rho_1}{2}\phi\coloneqq\rho_{\mathrm{avg}}+\rho_{\mathrm{dif}}\phi .
\end{equation}
\begin{remark}
	We have written the equation \eqref{eq:NS-CH_model_rho} in its more general variational formulation since $\JJ$ does not necessarily belong to $H^1(\Omega)^d$. It is clear from \eqref{rho} that $\rho_1\le \rho(\phi) \le \rho_2$ in $\Omega\times(0,T)$  is equivalent to $-1\le \phi \le 1$  in  $\Omega\times(0,T)$. Consequently, it is important the constraint  $\phi\in[-1,1]$ to preserve the physical meaning of the model because  the density of the mixture $\rho(\phi)$ must satisfy $\rho(\phi)\in[\rho_1,\rho_2]$.

\end{remark}

Finally, $\eta\in \mathcal{C}^0([-1,1])$ with $\eta(\phi)\ge C$ for certain $C>0$ and for all $\phi \in [-1,1]$ is the viscosity of the mixture, $\lambda>0$ is a constant related to the  energy density and $\varepsilon>0$ is a small parameter related to the thickness of the interface between the two fluids.

Since if $p$ is a pressure function solution of \eqref{NS-CH_model} then $p+C(t)$ is also solution for any time-dependent function $C(t)$, it is usual to consider the zero mean-value pressure constraint $\int_\Omega p=0$. 

We can consider the following variational formulation of problem \eqref{NS-CH_model}:  
Find $(\uu,p,\phi,\mu)$ such that $\uu\in L^\infty(0,T;L^2(\Omega)^d)\cap L^2(0,T;H_0^1(\Omega)^d)$, 
$p\in W^{-1,\infty}(0,T;L^2(\Omega))$ with $\int_\Omega p=0$, 
$\phi\in L^\infty(0,T;H^1(\Omega))$ 
with $-1\le \phi \le 1$ a.e. in  $\Omega\times(0,T)$, $\mu:\Omega\times(0,T)\to \mathbb{R}$ with  $\sqrt{M(\phi)} \nabla \mu \in L^2(0,T;L^2(\Omega))$, satisfying 
\begin{subequations} \label{eq:var_form_NS-CH}
	\begin{align}
	\label{eq:var_form_NS-CH_u}
	\langle\rho(\phi)\uu_t,\buu\rangle+\escalarLd{\left[(\rho(\phi)\uu-\rho_{\mathrm{dif}}M(\phi)\nabla\mu)\cdot\nabla\right]\uu}{\buu}&\nonumber\\+2\escalarLd{\eta(\phi)\DD\uu}{\DD\buu}-\escalarL{p}{\nabla\cdot\buu}-\escalarL{\mu}{\nabla\cdot(\phi\buu)}&=0,\\
	\label{eq:var_form_NS-CH_p}
	\escalarL{\nabla\cdot\uu}{\bp}&=0,\\
	\label{eq:var_form_NS-CH_phi}
	\langle\phi_t,\bphi\rangle+\escalarL{\nabla\cdot(\phi\uu)}{\bphi}+\escalarLd{M(\phi)\nabla \mu}{\nabla\bphi}&=0,\\
	\label{eq:var_form_NS-CH_mu}
	\lambda\varepsilon\escalarLd{\nabla\phi}{\nabla\bmu}+\frac{\lambda}{\varepsilon}\escalarL{f(\phi)}{\bmu}-\escalarL{\mu}{\bmu}&=0,
	\end{align}
\end{subequations}
for each $(\buu,\bp,\bmu,\bphi)\in (H_0^1(\Omega)\cap L^\infty(\Omega))^d\times L^2(\Omega)\times H^1(\Omega)\times H^1(\Omega)$. 
We have denoted $(f,g)=\int_\Omega f\,g$ as the $L^2(\Omega)$ scalar product and 
$$
\escalarLd{\eta(\phi)\DD\uu}{\DD\buu} = \int_\Omega\eta(\phi) \DD\uu :\DD\buu,
$$
where $:$ denotes the Frobenius inner product. 
\begin{proposition}
	The mass of the phase-field variable is conserved, because it holds
	$$
	\frac{d}{dt}\int_\Omega\phi=0.
	$$
	In particular, the mass of the mixture is conserved, because using \eqref{rho},
 	$$
	\int_\Omega \rho(\phi)
	= \vert\Omega\vert\rho_{\mathrm{avg}}+\rho_{\mathrm{dif}}\int_\Omega\phi
	= \vert\Omega\vert\rho_{\mathrm{avg}}+\rho_{\mathrm{dif}}\int_\Omega\phi_0
	=\int_\Omega \rho(\phi_0).
	$$
\end{proposition}
\begin{proof}
	Just test \eqref{eq:var_form_NS-CH_phi} by $\bphi=1$.
\end{proof}

\begin{proposition}
	Assuming a sufficiently regular solution of \eqref{eq:var_form_NS-CH_u}-\eqref{eq:var_form_NS-CH_mu}, the following energy law holds:
	\begin{equation}
		\label{energy_law_NS-CH}
	\frac{d}{dt} E(\uu,\phi)+2\int_\Omega\eta(\phi)\vert\DD\uu\vert^2+\int_\Omega M(\phi)\vert\nabla\mu\vert^2=0,
	\end{equation}
	where $\vert\DD\uu\vert^2=\sum_{i=1}^d\vert\DD\uu_i\vert^2$, with $\DD\uu_i$ denoting the $i$-th row of the stress tensor $\DD\uu$, and
	\begin{equation}
		\label{def:energy}
		E(\uu,\phi)\coloneqq\int_\Omega\rho(\phi)\frac{\vert\uu\vert^2}{2}+\frac{\lambda\varepsilon}{2}\int_\Omega\vert\nabla\phi\vert^2+\frac{\lambda}{\varepsilon}\int_\Omega F(\phi),
	\end{equation}
	where the first term is associated to the kinetic energy and the others to the potential energy. In particular, 
	the energy $E(\uu,\phi)$ is time decreasing because
	$$
	\frac{d}{dt} E(\uu,\phi)\le 0.
	$$
\end{proposition}
\begin{proof}
	We argue formally, by considering that all the functions that appear below are regular enough so that the expressions are true. Moreover, they are regarded as functions to be evaluated at $t\in(0,T)$, although, for clarity, we will omit it.
	
	If we test \eqref{eq:var_form_NS-CH_u}--\eqref{eq:var_form_NS-CH_mu} by $\buu=\uu$, $\bp=p$, $\bphi=\mu$ and $\bmu=\phi_t$ and we add up the expressions, we obtain:
	\begin{multline*}
	\escalarL{\rho(\phi)\uu_t}{\uu}+\lambda\varepsilon\escalarLd{\nabla\phi}{\nabla\phi_t}+\frac{\lambda}{\varepsilon}\escalarL{F'(\phi)}{\phi_t}\\+\escalarLd{\left[(\rho(\phi)\uu
	- \JJ)\cdot\nabla\right]\uu}{\uu}+2\int_\Omega\eta(\phi)\vert\DD\uu\vert^2 +\int_\Omega M(\phi)\vert\nabla\mu\vert^2=0.
	\end{multline*}

	Now, testing \eqref{eq:NS-CH_model_rho} by $\brho=\vert\uu\vert^2/2$, we have
	$$
	\escalarL{\partial_t\rho(\phi)}{\frac{\vert\uu\vert^2}{2}}-\escalarL{[(\rho(\phi)\uu-\JJ)\cdot\nabla]\uu}{\uu}=0.
	$$
	By adding the two previous expressions, the convective term $\escalarL{[(\rho(\phi)\uu-\JJ)\cdot\nabla]\uu}{\uu}$ cancels.
	Hence, taking into account that
	\begin{align*}
		\frac{d}{dt}\int_\Omega\rho(\phi)\frac{\vert\uu\vert^2}{2}&=\escalarL{\rho(\phi)\uu_t}{\uu}+\escalarL{\partial_t\rho(\phi)}{\frac{\vert\uu\vert^2}{2}},\\
		\frac12\frac{d}{dt}\int_\Omega\vert\nabla\phi\vert^2&=\escalarLd{\nabla\phi}{\nabla\phi_t},\\
		\frac{d}{dt}\int_\Omega F(\phi)&=\escalarL{F'(\phi)}{\phi_t},
	\end{align*}
	we can conclude that the energy law \eqref{energy_law_NS-CH} holds.
\end{proof}

\section{Structure-preserving scheme}
\label{sec:coupled_scheme}
In this section we develop a fully coupled discretization of the model \eqref{NS-CH_model} that preserves all properties at the discrete level, including the mass conservation, pointwise bounds of the phase-field and density of the mixture variables,  and the decreasing of the energy (also called energy-stability).

\subsection{Notation}
\label{sec:notation}
 We consider a finite element shape-regular triangular mesh $\T_h=\{K\}_{K\in \T_h}$ in the sense of Ciarlet, \cite{ciarlet2002finite}, of size $h$ over $\Omega$. We denote by $\E_h$ the set of the edges of $\T_h$ (faces if $d=3$) with $\Ehi$ the set of the \textit{interior edges} and $\Ehb$ the \textit{boundary edges}, i.e. $\E_h=\Ehi\cup\Ehb$.

 Now, we fix the following orientation over the mesh $\T_h$:
 \begin{itemize}
	\item For any interior edge $e\in\Ehi$ we set the associated unit normal vector $\nn_e$. In this sense, when referring to edge $e\in\Ehi$ we will denote by $K_e$ and $L_e$ the elements of $\T_h$ with $e=\partial K_e\cap\partial L_e$ and so that $\nn_e$ is exterior to $K_e$ pointing to $L_e$. 
	
	If there is no ambiguity, to abbreviate the notation  we will denote the previous elements $K_e$ and $L_e$ simply by $K$ and $L$, respectively, with the assumption that their naming is always with respect to the edge $e\in\Ehi$ and it may vary if we consider a different edge of $\Ehi$.
	\item For any boundary edge $e\in\Ehb$, the unit normal vector $\nn_e$ points outwards of the domain $\Omega$.
 \end{itemize}

 Therefore, we can define the \textit{average} $\media{\cdot}$ and the \textit{jump} $\salto{\cdot}$ of a function $v$ on an edge $e\in\E_h$ as follows:
\begin{equation*}
		\media{v}\coloneqq
		\begin{cases}
			\dfrac{\vK+\vL}{2}&\text{if } e\in\Ehi,\,  e=K\cap L,\\
			\vK&\text{if }e\in\Ehb,\, e\subset K,
		\end{cases}
		\qquad
		\salto{v}\coloneqq
		\begin{cases}
			\vK-\vL&\text{if } e\in\Ehi,\, e=K\cap L,\\
			\vK&\text{if }e\in\Ehb,\, e\subset K.
		\end{cases}
\end{equation*}

We denote by $\Pd_k(\T_h)$ and $\Pc_k(\T_h)$ the spaces of finite element discontinuous and continuous functions, respectively, which are polynomials of degree $k\ge0$ when restricted to the elements $K$ of $\T_h$.
In this sense, we will denote the broken differential operators (see \cite{riviere_discontinuous_2008,di_pietro_mathematical_2012}) the same way as the standard differential operators in the absence of ambiguity.

Moreover, we take an equispaced partition $0=t_0<t_1<\cdots<t_N=T$ of the time domain $[0,T]$ with $\Delta t=t_{m+1}-t_m$ the time step. Also, for any function $v$ depending on time, we denote 
$v^{m+1}\simeq v(t_{m+1})$ and the discrete time derivative operator $v_t(t_{m+1})\simeq\delta_t v^{m+1}:=(v^{m+1}-v^m)/\Delta t$.

Finally, we set the following notation for the positive and negative parts of a function $v$:
	$$
	v_\oplus\coloneqq\frac{|v|+v}{2}=\max\{v,0\},
	\quad
	v_\ominus\coloneqq\frac{|v|-v}{2}=-\min\{v,0\},
	\quad
	v=v_\oplus - v_\ominus.
	$$

\subsection{Discrete scheme}

Following the ideas of \cite{acosta-soba_upwind_2022,acosta-soba_KS_2022,acosta2023structure} we define the projections $\Pi_0\colon L^1(\Omega) \longrightarrow \Pd_0(\T_h)$, $\Pi_1\colon L^1(\Omega)\longrightarrow \Pc_1(\T_h)$ and 
$\Pih_1\colon L^1(\Omega)\longrightarrow \Pc_1(\T_h)$ 
as follows:
\begin{align}
	\label{eq:esquema_DG_Pi0}
	\escalarL{\Pi_0 g}{\overline{w}}
	&=\escalarL{g}{\overline{w}}, 
	&\forall\,\overline{w}\in \Pd_0(\T_h),
	\\
	\label{eq:esquema_DG_Pi1}
	\escalarL{\Pi_1 g}{\overline{v}}
	&=
	\escalarL{g}{\overline{v}},&\forall\,\overline{v}\in \Pc_1(\T_h),
	\\
	\label{eq:esquema_DG_Pih1}
	\escalarML{\Pih_1 g}{\overline{v}}
&=\escalarL{g}{\overline{v}},&\forall\,\overline{v}\in \Pc_1(\T_h),
\end{align}
where $\escalarL{\cdot}{\cdot}$ denotes the usual scalar product in $L^2(\Omega)$. In addition, $\escalarML{\cdot}{\cdot}$ denotes the mass-lumping scalar product in $\Pc_1(\T_h)$ resulting from using the trapezoidal rule to approximate the scalar product in $L^2(\Omega)$ (see, for instance, \cite{quarteroni2008numerical}). Therefore, for any elements $\varphi, \psi\in\Pc_1(\T_h)$ this scalar product can be defined as $$\escalarML{\varphi}{\psi}=\frac{1}{3}\sum_{K\in\T_h}\vert K\vert \sum_{j=1}^3\varphi(x_{j,K})\psi(x_{j,K}),$$ where $x_{j,K}$ are the nodes of the element $K\in\T_h$ for every $j\in\{1,2,3\}$. These projections \eqref{eq:esquema_DG_Pi0}--\eqref{eq:esquema_DG_Pih1} are well defined for every function $g\in L^1(\Omega)$, notice that $\overline{\varphi}\in\Pd_k(\T_h)$ imply $\overline{\varphi}|_K\in\mathcal{C}^\infty(K)$ for every $K\in\T_h$ and, therefore, $\escalarL{g}{\overline{\varphi}}=\sum_{K\in\T_h}\int_K g\overline{\varphi} < \infty$.

We propose the following numerical scheme: find $\uu^{m+1}\in \Uh$, $p^{m+1}\in\Ph$ 
with $\int_\Omega p^{m+1}=0$, $\phi^{m+1}\in \Pd_0(\T_h)$ and $\mu^{m+1}\in \Pc_1(\T_h)$ such that
\begin{subequations}
	\label{esquema_DG_NS-CH}
	\begin{align}
		\label{eq:esquema_DG_NS-CH_u}
		\escalarL{\rho(\Pih_1\phi^{m})\delta_t\uu^{m+1}}{\buu}
		+\escalarL{\left[\left(\rho(\Pih_1\phi^{m})u^m-\JJ^m_h\right)\cdot\nabla\right]\uu^{m+1}}{\buu}
		&
		\nonumber\\
		+2\escalarLd{\eta(\phi^{m})\DD\uu^{m+1}}{\DD\buu}
		-\escalarL{p^{m+1}}{\nabla\cdot\buu}
		+\ch{\phi^{m+1}}{\Pi_0\mu^{m+1}}{\buu}
		&
		\nonumber\\
		+\th{\uu^{m+1}}{\uu^{m}}{\Pih_1\phi^{m+1}}{\Pih_1\phi^m}{\mu^m}{\buu}
		+\sh{\uu^{m+1}}{\phi^{m+1}}{\Pi_0\mu^{m+1}}{\buu}&=0,
		\\
		\label{eq:esquema_DG_NS-CH_p}
		\escalarL{\nabla\cdot\uu^{m+1}}{\bp}&=0,
		\\
		\label{eq:esquema_DG_NS-CH_phi}
		\escalarL{\delta_t\phi^{m+1}}{\bphi}
		+\aupw{\uu^{m+1}}{\phi^{m+1}}{\bphi}
		+\bupw{\nabla_{\nn}^0\mu^{m+1}}{M(\phi^{m+1})}{\bphi}&=0,\\
		\label{eq:esquema_DG_NS-CH_mu}
		\lambda\varepsilon\escalarLd{\nabla (\Pih_1 \phi^{m+1})}{\nabla\bmu}
		+\frac{\lambda}{\varepsilon}\escalarL{f(\Pih_1\phi^{m+1}, \Pih_1\phi^m)}{\bmu}-\escalarML{\mu^{m+1}}{\bmu}&=0,
	\end{align}
\end{subequations}
for each $\buu\in \Uh$, $\bp\in\Ph$, $\bphi\in\Pd_0(\T_h)$, $\bmu\in\Pc_1(\T_h)$,
where 
$$\JJ^m_h=\rho_\mathrm{dif}M(\Pih_1\phi^m)\Pi_1(\nabla\mu^m),$$ 
and
\begin{equation}
\label{def:convex_splitting}
f(\phi_1,\phi_0)\coloneqq F_i'(\phi_1) + F_e'(\phi_0) \text{ with } F_i(\phi)\coloneqq \phi^2+\frac{1}{4},\,F_e(\phi)\coloneqq \frac{1}{4}\phi^4-\frac{3}{2}\phi^2
\end{equation}
such that $F(\phi)=F_i(\phi)+F_e(\phi)$ is a convex splitting discretization of the Ginzburg-Landau double well potential $F(\phi)$ for any  $\phi\in [-1,1]$.

Also, $(\Uh,\Ph)$ is a compatible ``inf-sup" pair of finite-dimensional spaces satisfying that $\Uh\subset(\mathcal{C}^0(\overline\Omega)\cap H^1_0(\Omega))^d$ and $\Pd_0(\T_h) \subset \Ph $.
In fact, the restriction $\Pd_0(\T_h)\subset\Ph$ is needed in order to guarantee the local incompressibility of $\uu^{m+1}$ in the following sense:
\begin{equation}
	\label{local_incompressibility}
	\sum_{e\in\Ehi}\int_{e} (\uu^{m+1}\cdot\nn_e) \salto{\bp}=0,\quad\forall \,\bp\in\Pd_0(\T_h),
\end{equation}
which can be derived integrating by parts in \eqref{eq:esquema_DG_NS-CH_p}.
This constraint will allow us to preserve the pointwise bounds of $\phi^{m+1}$, see Theorem \ref{thm:discrete_maximum_principle} below. Notice that the discretization of the pressure and the divergence term~(\ref{eq:esquema_DG_NS-CH_p}) is the standard Stokes DG approach~\cite{riviere_discontinuous_2008, di_pietro_mathematical_2012} for continuous velocity and discontinuous pressure.

\begin{remark}
	\label{rmk:inf-sup_spaces}
Some possible choices  of compatible spaces  $(\Uh,\Ph)$ are the following (see \cite{boffi2013mixed,ern_theory_2010} for the details):
\begin{itemize}
	\item $(\Uh,\Ph)=((\Pc_2(\T_h)\cap H_0^1(\Omega))^d,\Pd_0(\T_h))$ which is stable for $d=2$ but not for $d=3$.
	\item $(\Uh,\Ph)=((\Pb_2(\T_h)\cap H^1_0(\Omega))^d,\Pd_1(\T_h))$ which is stable for $d=2,3$ but requires a higher computational effort. Here, $\Pb_2(\T_h)$ denotes the $\Pc_2(\T_h)$ space enriched with a bubble by elements of order 3.
	\item $(\Uh,\Ph)=((\mathcal{Q}_2(\T_h)\cap H^1_0(\Omega))^d,\Pd_1(\T_h))$. Here, $\mathcal{Q}_2(\T_h)$ denotes the standard quadrilateral finite element space of order 2.
	
\end{itemize}
\end{remark}

Notice that, for any choice of this pair $(\Uh,\Ph)$, the error bounds are expected to be determined by the lowest accuracy approximation
of the phase-field function by $\Pd_0(\T_h)$.

Moreover,
$\ch{\phi}{\mu}{\buu}$ is a centered discretization of the term $\escalarL{\phi\nabla\mu}{\buu}=-\escalarL{\mu}{\nabla\cdot(\phi\buu)}$ in \eqref{eq:var_form_NS-CH_u}
defined as
\begin{equation}
	\label{centered_discretization}
	\ch{\phi}{\mu}{\buu}\coloneqq 
	- \int_\Omega\nabla\cdot(\phi\buu) \mu
	-\sum_{e\in\Ehi}\int_e (\buu\cdot\nn_e)\media{\phi}\salto{\mu},
\end{equation}
where 
the second term is a consistent stabilization term depending on the jumps of $\mu$ on the interior edges of the mesh $\T_h$.

In \eqref{eq:esquema_DG_NS-CH_phi} we have considered two different upwind formulas, the classical upwind 
\begin{align}
	\label{def:aupw}
	\aupw{\uu}{\phi}{\bphi}&\coloneqq \sum_{e\in\Ehi, e=K\cap L}\int_e\left( (\uu\cdot\nn_e)_\oplus\phi_K - (\uu\cdot\nn_e)_\ominus\phi_L\right)\salto{\bphi}
\end{align}
whose properties were discussed in \cite{acosta-soba_upwind_2022}, and 
$$\bupw{\nabla_{\nn}^0\mu }{M(\phi)}{\bphi},$$
 which follows the ideas introduced in \cite{acosta-soba_KS_2022,acosta2023structure}, and which will be detailed in the Subsection \ref{sec:def_bh}. 

Finally, we have introduced in \eqref{eq:esquema_DG_NS-CH_u} two consistent stabilizations terms:
\begin{equation}
	\label{def:sh1}
	\th{\uu_1}{\uu_0}{\phi_1}{\phi_0}{\mu}{\buu}\coloneqq 
	\frac{1}{2}\Big\{\escalarL{\delta_t\rho(\phi_1)}{\uu_1\cdot \buu}-\escalarLd{\rho(\phi_0)\uu_0-\rho_\mathrm{dif}M(\phi_0)\Pi_1(\nabla\mu)}{\nabla(\uu_1\cdot \buu)}\Big\},
\end{equation}
which, following the ideas of \cite{tierra_guillen_abels_2014}, can be interpreted as a residual to the equation \eqref{eq:NS-CH_model_rho}; and
\begin{equation}
\label{def:sh2}
	\sh{\uu}{\phi}{\mu}{\buu}\coloneqq
	-\frac{1}{2}\sum_{e\in\Ehi}\int_e(\buu\cdot\nn_e)\, \sign{\uu\cdot\nn_e}\salto{\phi}\salto{\mu},
\end{equation}
which 
is introduced to control the influence of the upwind term $\aupw{\uu^{m+1}}{\phi^{m+1}}{\bphi}$ in \eqref{eq:esquema_DG_NS-CH_phi}. This latter stabilization together with the centered approximation 
 $\ch{\phi^{m+1}}{\Pi_0\mu^{m+1}}{\buu}$
 of the phase-field force in the momentum equation \eqref{eq:esquema_DG_NS-CH_u}, cancel the effect of the transport of the phase-field function by the mean velocity $\uu^{m+1}$ and allow us to obtain a discrete energy inequality, see Lemma \ref{lemma:discrete_energy} below.

To start the algorithm we take $\phi^0=\Pi_0\phi_0$ where $\phi_0$ is the continuous initial data, which  satisfies $\phi_0\in[-1,1]$. Notice that, one also has $\phi^0\in[-1,1]$.

\begin{remark}
	Observe that the $0$-mean value constraint on the pressure has been removed from the discrete formulation \eqref{esquema_DG_NS-CH}. This constraint will be enforced
	in practice by using an additional penalty term, see Section~\ref{sec:numerical_experiments} below.
\end{remark}

\subsubsection{Definition of the  upwind bilinear form $\boldmath\bupw{\cdot}{\cdot}{\cdot}$}
\label{sec:def_bh}
In order to define the upwind bilinear form $\bupw{\cdot}{\cdot}{\cdot}$ we follow the ideas of \cite{acosta-soba_KS_2022,acosta2023structure}.

First, we split the mobility function $M(z)$ for $z\in\R$ into its increasing and decreasing parts, denoted respectively by $M^\uparrow(z)$ and $M^\downarrow(z)$, as follows:
\begin{subequations}
	\begin{align*}
		M^\uparrow(z)&
		=\int_{-1}^{\min(z,1)}M'(s)_\oplus ds=\int_{-1}^{\min(z,1)}(-2s)_\oplus ds,
		\\
		M^\downarrow(z)&
		=-\int_{-1}^{\min(z,1)} M'(s)_\ominus ds=-\int_{-1}^{\min(z,1)} (-2s)_\ominus ds
	\end{align*}
\end{subequations}
Therefore,
\begin{align}
	\label{movilidad_partes_creciente_decreciente}
	M^\uparrow(z)=
	\begin{cases}
		M(z) & \text{if }z\le 0\\[0.2em]
		M(0) & \text{if } z>0
	\end{cases},\quad
	M^\downarrow(z)=
	\begin{cases}
		0 & \text{if }z\le 0\\
		M(z)-M(0) & \text{if } z>0
    \end{cases}.
\end{align}
Notice that $M^\uparrow(z)+ M^\downarrow(z) = M(z)$.

Following the work in \cite{acosta-soba_upwind_2022}, we can define the following upwind form for any  $\phi, \bphi\in\Pd_0(\T_h)$ and $\mu\in\Pc_1(\T_h)$:
\begin{align}
	\label{def:bupw_general}
&\bupw{-\nabla_{\nn}\mu}{M(\phi)}{\bphi}\coloneqq\nonumber\\&\sum_{e\in\Ehi,e=K\cap L}\int_e\left((-\nabla_{\nn}\mu)_\oplus(M^\uparrow(\phi_K)+M^\downarrow(\phi_L))_\oplus-(-\nabla_{\nn}\mu)_\ominus(M^\uparrow(\phi_L)+M^\downarrow(\phi_K))_\oplus\right)\salto{\bphi},
\end{align}
where $\nabla_{\nn}\mu\coloneqq
\media{\nabla\mu}\cdot\nn_e$ on every $e\in\E_h$.

Nonetheless, if we want to ensure a discrete energy law, as was done in \cite{acosta-soba_KS_2022,acosta2023structure}, we need to introduce the following hypothesis:
\begin{hyp}
	\label{hyp:mesh_n}
	The mesh $\T_h$ of $\overline\Omega$ is structured in the sense that, for any interior interface $e=K\cap L\in\Ehi$,  the line between the barycenters of $K$ and $L$ is orthogonal to $e$.
\end{hyp}

Under this hypothesis, we can consider the following consistent approximation on every $e\in\E_h^i$, as done in \cite{acosta-soba_KS_2022,acosta2023structure}:
\begin{equation}
\label{eq:approx_gradn}
\nabla\mu\cdot\nn_e \simeq\frac{-\salto{\Pi_0\mu}}{\mathcal{D}_e(\T_h)}\coloneqq \nabla_{\nn}^0\mu|_{e},
\end{equation}
where $\mathcal{D}_e(\T_h)$ is the distance between the barycenters of the triangles of the mesh $\T_h$ that share $e\in\Ehi$.

Therefore, we can extend the definition of the upwind bilinear form \eqref{def:bupw_general} as follows:
\begin{align}
	&\bupw{-\nabla_{\nn}^0\mu}{M(\phi)}{\bphi} \nonumber
	\\& 
	\label{def:bupw_saltos}
=	\sum_{e\in\Ehi,e=K\cap L}\frac{1}{\mathcal{D}_e(\T_h)}\int_e\left((\salto{\Pi_0\mu})_\oplus(M^\uparrow(\phi_K)+M^\downarrow(\phi_L))_\oplus-(\salto{\Pi_0\mu})_\ominus(M^\uparrow(\phi_L)+M^\downarrow(\phi_K))_\oplus\right)\salto{\bphi}.
\end{align}
This upwind approximation allows us to obtain both a discrete maximum principle and an energy-stability property as shown in \cite{acosta2023structure} for a tumor model based on the Cahn-Hilliard equation with degenerate mobility.

\begin{remark}
	Notice that the upwind bilinear form $\aupw{\uu}{\phi}{\bphi}$ given in \eqref{def:aupw}, can be seen as a particular case of $\bupw{\cdot}{\cdot}{\cdot}$ given in \eqref{def:bupw_general}, changing $M(\phi)$ by $\phi$,  but now we have not truncated the transported variable $\phi$.
	In fact, it is not necessary to truncate  $\phi$ in $\aupw{\uu}{\phi}{\bphi}$ to preserve the pointwise bounds of $\phi$ due to the local incompressibility of $\uu$ (see \cite{acosta-soba_upwind_2022} for a more detailed explanation).
\end{remark}

\subsubsection{Properties of the scheme \eqref{esquema_DG_NS-CH}}
\label{sec:properties_scheme}
\begin{proposition}[Mass conservation]
	The mass of the phase-field variable and its regularization are conserved. In fact, one has 
	$$
	\int_\Omega\phi^{m+1}=\int_\Omega \phi^m,\quad \int_\Omega\Pih_1\phi^{m+1}=\int_\Omega \Pih_1\phi^m.
	$$

	As a consequence, since $\rho(\phi)$ is linear with respect to $\phi$, the mass of the mixture is also conserved, 
	$$
   \int_\Omega \rho(\phi^{m+1})=\int_\Omega \rho(\phi^{m}),\quad \int_\Omega \rho(\Pih_1\phi^{m+1})=\int_\Omega \rho(\Pih_1\phi^{m}).
   $$
\end{proposition}
\begin{proof}
	Just need to take $\bphi=1$ in \eqref{eq:esquema_DG_NS-CH_phi} and consider the definitions of the regularization $\Pih_1$ given in \eqref{eq:esquema_DG_Pih1}, and the density of the mixture $\rho(\phi)$ given in \eqref{rho}.
\end{proof}

\begin{theorem}[Pointwise bounds of the phase-field variable]
	\label{thm:discrete_maximum_principle}
	Provided that $\phi^m\in[-1,1]$ in $\Omega$, any solution $\phi^{m+1}$ of \eqref{esquema_DG_NS-CH}  
   and 
  $\Pih_1\phi^{m+1}$ satisfy:
  $\phi^{m+1}, \Pih_1\phi^{m+1}\in[-1,1]$ in $\Omega$.
\end{theorem}
\begin{proof}
	To prove that $\phi^{m+1}\ge -1$ in $\Omega$ we may take the following $\Pd_0(\T_h)$ test function
	\begin{align*}
		\bphi^*=
		\begin{cases}
			(\phi^{m+1}_{K^*}+1)_\ominus &\text{in }K^*\\
			0&\text{out of }K^*
		\end{cases},
	\end{align*}
	where $K^*$ is an element of $\T_h$ such that  $\phi_{K^*}^{m+1}=\min_{K\in\T_h}\phi_{K}^{m+1}$. We denote $\nn_{K^*}$ the normal vector exterior to $K^*$. Then,  since $\phi_L^{m+1}\ge \phi_{K^*}^{m+1}$ we can assure, using the local incompressibility constraint \eqref{local_incompressibility}, that
	\begin{align*}
	\aupw{\uu^{m+1}}{\phi^{m+1}}{\bphi^*}&
	= \sum_{e\in\Ehi}\int_e\left( (\uu^{m+1}\cdot\nn_e)_\oplus\phi_{K}^{m+1} 
	- (\uu^{m+1}\cdot\nn_e)_\ominus\phi_L^{m+1}\right) \salto{\bphi^*}
	\\&
	=  \sum_{e\in\Ehi, e=K^*\cap L}\int_e\left( (\uu^{m+1}\cdot\nn_{K^*})_\oplus\phi_{K^*}^{m+1} 
	- (\uu^{m+1}\cdot\nn_{K^*})_\ominus\phi_L^{m+1}\right)(\phi_{K^*}^{m+1}+1)_\ominus
	\\
	&\le
	\sum_{e\in\E_h^i,e\subset K^*}\int_e(\uu^{m+1}\cdot\nn_{K^*})\phi_{K^*}^{m+1}(\phi_{K^*}^{m+1}+1)_\ominus\\&=\sum_{e\in\E_h^i}\int_e(\uu^{m+1}\cdot\nn_e)\salto{\phi^{m+1}\bphi^*}=0.
	\end{align*}
	
	On the other hand, using that the positive part is an increasing function and that
	$$M^\uparrow(\phi_L^{m+1}) \ge M^\uparrow(\phi_{K^*}^{m+1})\quad
	\hbox{and} \quad
	M^\downarrow(\phi_L^{m+1}) \le M^\downarrow(\phi_{K^*}^{m+1}),
	$$
	we can obtain (see \cite{acosta-soba_upwind_2022,acosta2023structure})
	$$
	\bupw{\nabla_{\nn}^0\mu^{m+1}}{M(\phi^{m+1})}{\bphi^*}\le 0.
	$$
	
	Consequently, $\vert K^*\vert\delta_t \uKs^{m+1}(\uKs^{m+1}+1)_\ominus \ge0$. Therefore,
	$$
	0
	\le |K^*|(\delta_t (\phi_{K^*}^{m+1}+1))(\phi_{K^*}^{m+1}+1)_\ominus 
	=
	-\frac{|K^*|}{\Delta t}\left((\phi_{K^*}^{m+1}+1)_\ominus ^2+(\phi_{K^*}^m+1)(\phi_{K^*}^{m+1}+1)_\ominus \right)
	\le 0,
	$$
	which implies, since $\phi_{K^*}^m\ge -1$, that $(\phi_{K^*}^{m+1}+1)_\ominus =0$. Hence, $\phi^{m+1}\ge-1$
	 in $\Omega$.

	Similarly, taking the following $\Pd_0(\T_h)$ test function
	\begin{align*}
		\bphi^*=
		\begin{cases}
			(\phi^{m+1}_{K*}-1)_\oplus &\text{in }K^*\\
			0&\text{out of }K^*
		\end{cases},
	\end{align*}
	where $K^*$ is an element of $\T_h$ such that  $\phi_{K^*}^{m+1}=\max_{K\in\T_h}\phi_{K}^{m+1}$, we can arrive at $\phi^{m+1}\le 1$ in $\Omega$.

	Finally, $\Pih_1\phi^{m+1}\in[-1,1]$ in $\Omega$ is a direct consequence of the definition of the projection  $\Pih_1$ given in \eqref{eq:esquema_DG_Pih1}.
\end{proof}

The next Corollary is a direct consequence of the previous result.
\begin{corollary}[Pointwise bounds of the mixture density]
	\label{cor:discrete_maximum_principle}
	Provided that $\rho(\phi^m)\in[\rho_1,\rho_2]$ in $\Omega$, the density of the mixture satisfies $\rho(\phi^{m+1}),\rho(\Pih_1\phi^{m+1})\in[\rho_1,\rho_2]$ in $\Omega$.
\end{corollary}

The following Lemma is a technical result that we are going to use when computing the discrete energy law.
\begin{lemma}
	\label{lemma:discrete_energy}
	The following expression holds
	\begin{equation}
		\label{upwind_stabilization_energy}
		\aupw{\uu^ {m+1}}{\phi^{m+1}}{\Pi_0\mu^{m+1}}
		+\ch{\phi^{m+1}}{\Pi_0\mu^{m+1}}{\uu^{m+1}}
		+\sh{\uu^{m+1}}{\phi^{m+1}}{\Pi_0\mu^{m+1}}{\uu^{m+1}}=0.
	\end{equation}
\end{lemma}
\begin{proof}
	First, notice that we can rewrite the term $\aupw{\uu^ {m+1}}{\phi^{m+1}}{\Pi_0\mu^{m+1}}$ as follows
	\begin{multline*}
		\aupw{\uu^ {m+1}}{\phi^{m+1}}{\Pi_0\mu^{m+1}}=\sum_{e\in\E_h}\int_e(\uu^{m+1}\cdot\nn_e)\media{\phi^{m+1}}\salto{\Pi_0\mu^{m+1}}\\+\frac{1}{2}\sum_{e\in\Ehi}\int_e\vert \uu^{m+1}\cdot\nn_e\vert\salto{\phi^{m+1}}\salto{\Pi_0\mu^{m+1}}.
	\end{multline*}
	Then, by definition and due to $\phi^{m+1}\in\Pd_0(\T_h)$,
	\begin{align*}
		\ch{\phi^{m+1}}{\Pi_0\mu^{m+1}}{\uu^{m+1}}&
		=-\int_\Omega(\nabla\cdot\uu^{m+1})\phi^{m+1}\Pi_0\mu^{m+1}\\ &\quad-\sum_{e\in\E_h}\int_e (\uu^{m+1}\cdot\nn_e)\media{\phi^{m+1}}\salto{\Pi_0\mu^{m+1}},\\
		\sh{\uu^{m+1}}{\phi^{m+1}}{\Pi_0\mu^{m+1}}{\uu^{m+1}}&=-\frac{1}{2}\sum_{e\in\Ehi}\int_e\vert\uu^{m+1}\cdot\nn_e\vert\salto{\phi^{m+1}}\salto{\Pi_0\mu^{m+1}}.
	\end{align*}
	Finally, using \eqref{eq:esquema_DG_NS-CH_p},
	$$
	\ch{\phi^{m+1}}{\Pi_0\mu^{m+1}}{\uu^{m+1}}=-\sum_{e\in\E_h}\int_e (\uu^{m+1}\cdot\nn_e)\media{\phi^{m+1}}\salto{\Pi_0\mu^{m+1}},
	$$
	what yields \eqref{upwind_stabilization_energy}.
\end{proof}

\begin{theorem}[Discrete energy law]
	\label{thm:discrete_energy_law_NS-CH}
	The following discrete energy law holds:
	\begin{align}
	\label{discrete_energy_law_NS-CH}
	\delta_tE(\uu^{m+1},\Pih_1\phi^{m+1})&+2\escalarLd{\eta(\phi^{m+1})\DD\uu^{m+1}}{\DD\uu^{m+1}}
	+\bupw{-\nabla_{\nn}^0\mu^{m+1}}{M(\phi^{m+1})}{\Pi_0\mu^{m+1}}\nonumber\\
	&+\frac{\Delta t}{2}\int_\Omega\rho(\Pih_1\phi^m)\vert \delta_t \uu^{m+1}\vert^2+\frac{\Delta t\lambda\varepsilon}{2}\int_\Omega\vert\delta_t\nabla\Pih_1\phi^{m+1}\vert^2\nonumber\\&+\frac{\lambda}{\varepsilon}\int_\Omega\left(f(\Pih_1\phi^{m+1}, \Pih_1\phi^m)\delta_t\Pih_1\phi^{m+1}-F(\Pih_1\phi^{m+1})\right)
	=0,
	\end{align}
	where the energy functional $E(\uu,\phi)$ is defined in \eqref{def:energy}.
\end{theorem}
\begin{proof}
	First, take $\buu=\uu^{m+1}$ and $\bp=p^{m+1}$ in \eqref{eq:esquema_DG_NS-CH_u}--\eqref{eq:esquema_DG_NS-CH_p}. Consider that
	\begin{align}
		\label{proof:discrete_energy_law_NS-CH_1}
		\escalarL{\rho(\Pih_1\phi^{m})\delta_t\uu^{m+1}}{\uu^{m+1}}&=\frac{1}{2}\int_\Omega\rho(\Pih_1\phi^m)\delta_t \vert \uu^{m+1}\vert^2+\frac{\Delta t}{2}\int_\Omega\rho(\Pih_1\phi^m)\vert \delta_t \uu^{m+1}\vert^2,
	\end{align}
	and, by definition of $\th{\cdot}{\cdot}{\cdot}{\cdot}{\cdot}{\cdot}$ given in \eqref{def:sh1},
	\begin{align}
		\label{proof:discrete_energy_law_NS-CH_2}
		\frac{1}{2}\int_\Omega\delta_t\left(\rho(\Pih_1\phi^{m+1})\right)\vert\uu^{m+1}\vert^2&=\escalarL{\left[\left(\rho(\Pih_1\phi^{m})u^m-\JJ^m_h\right)\cdot\nabla\right]\uu^{m+1}}{\uu^{m+1}}\nonumber\\&\quad+\th{\uu^{m+1}}{\uu^{m}}{\Pih_1\phi^{m+1}}{\Pih_1\phi^m}{\mu^m}{\uu^{m+1}}.
	\end{align}
	Then, using \eqref{proof:discrete_energy_law_NS-CH_1} and \eqref{proof:discrete_energy_law_NS-CH_2} we can arrive at the following expression
	\begin{multline}
	 \label{eq:energia_DG_u}
	\delta_t\int_\Omega\rho(\Pih_1\phi^{m+1})\frac{\vert\uu^{m+1}\vert^2}{2}
	+\frac{\Delta t}{2}\int_\Omega\rho(\Pih_1\phi^m)\vert \delta_t \uu^{m+1}\vert^2+2\escalarLd{\eta(\phi^{m+1})\DD\uu^{m+1}}{\DD\uu^{m+1}}\\
	+\ch{\phi^{m+1}}{\Pi_0\mu^{m+1}}{\uu^{m+1}}+\sh{\uu^{m+1}}{\phi^{m+1}}{\Pi_0\mu^{m+1}}{\uu^{m+1}}=0.
	\end{multline}

	Now, if we test \eqref{eq:esquema_DG_NS-CH_phi}--\eqref{eq:esquema_DG_NS-CH_mu} with $\bphi=\Pi_0\mu^{m+1}$ and $\bmu=\delta_t\Pih_1\phi^{m+1}$ and we add the resulting expressions and \eqref{eq:energia_DG_u}, we obtain, using \eqref{upwind_stabilization_energy},
	\begin{multline*}
		\delta_t\int_\Omega\rho(\Pih_1\phi^{m+1})\frac{\vert\uu^{m+1}\vert^2}{2}
		+\frac{\Delta t}{2}\int_\Omega\rho(\Pih_1\phi^m)\vert \delta_t \uu^{m+1}\vert^2+2\escalarLd{\eta(\phi^{m+1})\DD\uu^{m+1}}{\DD\uu^{m+1}}\\+\escalarL{\delta_t\phi^{m+1}}{\Pi_0\mu^{m+1}}+\bupw{-\nabla_{\nn}^0\mu^{m+1}}{M(\phi^{m+1})}{\Pi_0\mu^{m+1}}
		+\lambda\varepsilon\escalarLd{\nabla \Pih_1 \phi^{m+1}}{\delta_t\nabla\Pih_1\phi^{m+1}}\\+\frac{\lambda}{\varepsilon}\escalarL{f(\Pih_1\phi^{m+1},\Pih_1 \phi^m)}{\delta_t\Pih_1\phi^{m+1}}-\escalarML{\mu^{m+1}}{\delta_t\Pih_1\phi^{m+1}}=0.
	\end{multline*}

	Finally, the following equalities
	\begin{align*}
		\escalarL{\delta_t\phi^{m+1}}{\Pi_0\mu^{m+1}}&=\escalarL{\delta_t\phi^{m+1}}{\mu^{m+1}}=\escalarML{\delta_t\Pih_1\phi^{m+1}}{\mu^{m+1}},\\
		\lambda\varepsilon\escalarLd{\nabla \Pih_1 \phi^{m+1}}{\delta_t\nabla\Pih_1\phi^{m+1}}&=\frac{\lambda\varepsilon}{2}\delta_t\int_\Omega\vert\nabla\Pih_1\phi^{m+1}\vert^2+\frac{\Delta t\lambda\varepsilon}{2}\int_\Omega\vert\delta_t\nabla\Pih_1\phi^{m+1}\vert^2,
	\end{align*}
	yield \eqref{discrete_energy_law_NS-CH}.
\end{proof}

Using the definition of the upwind form $\bupw{\cdot} {\cdot}{\cdot}$ and the standard procedure for the convex-splitting technique (see e.g. \cite{eyre_1998_unconditionally,guillen-gonzalez_linear_2013}), one can show the following Lemma.
\begin{lemma}
	\label{lemma:discrete_energy_inequalities}
	The following two inequalities hold:
	\begin{align}
		\bupw{-\nabla_{\nn}^0\mu^{m+1}}{M(\phi^{m+1})}{\Pi_0\mu^{m+1}}&\ge 0,\\
		\int_\Omega\left(f(\Pih_1\phi^{m+1},\Pih_1 \phi^m)\delta_t\Pih_1\phi^{m+1}-\delta_t F(\Pih_1\phi^{m+1})\right)&\ge 0.
	\end{align}
\end{lemma}

The following result is a direct consequence of Theorem \ref{thm:discrete_energy_law_NS-CH} and Lemma \ref{lemma:discrete_energy_inequalities}.
\begin{corollary}[Discrete energy stability]
	\label{cor:discrete_energy_law_NS-CH}
	The scheme \eqref{esquema_DG_NS-CH} satisfies
	\begin{equation}
		\delta_tE(\uu^{m+1},\Pih_1\phi^{m+1})+2\escalarLd{\eta(\phi^{m+1})\DD\uu^{m+1}}{\DD\uu^{m+1}}+\bupw{-\nabla_{\nn}^0\mu^{m+1}}{M(\phi^{m+1})}{\Pi_0\mu^{m+1}}\le 0.
	\end{equation}
	In particular, scheme \eqref{esquema_DG_NS-CH} is unconditionally energy stable, i.e., $\delta_tE(\uu^{m+1},\Pih_1\phi^{m+1})\le 0$.
\end{corollary}

	The scheme \eqref{esquema_DG_NS-CH} is nonlinear so we will need to approximate its solution by means of an iterative procedure such as the nonsmooth Newton's method (see \cite{clarke1990optimization}).
	
	However, the function $\sign{\phi}$ that appears in the stabilization term $\sh{\cdot}{\cdot}{\cdot}{\cdot}$ is not subdifferentiable at $\phi=0$ and, although it is rare in practice that $\phi=0$ holds exactly due to round-off errors, one might eventually find convergence issues. In this case, several approaches can be carried out to improve the convergence of the algorithm. For instance, one may use an iterative procedure that does not rely on the Jacobian of the whole system such as a fixed point algorithm. Conversely, if we want to use a higher order procedure depending on the Jacobian like the nonsmooth Newton's method, one may avoid the use of the $\sign$ function regularizing the term $\sh{\cdot}{\cdot}{\cdot}{\cdot}$ as follows
	\begin{equation}
		\label{def:regularization_sh}
		\shd{\uu}{\phi}{\mu}{\buu}\coloneqq\frac{1}{2}\sum_{e\in\Ehi}\int_e(\buu\cdot\nn_e)\frac{\uu\cdot\nn_e}{\vert\uu\cdot\nn_e\vert+\delta}\salto{\Pi_0\mu}\salto{\phi},
	\end{equation}
	for $\delta>0$ small. This modification preserves the mass conservation and the pointwise bounds but introduces a modification in the discrete energy law, see Theorem \ref{thm:regularized_energy}.

The following result can be proved using the same procedure in  Theorem~\ref{thm:discrete_energy_law_NS-CH} and Corollary~\ref{cor:discrete_energy_law_NS-CH}.
\begin{theorem}
	\label{thm:regularized_energy}
	If we regularize the stabilization term $\sh{\cdot}{\cdot}{\cdot}{\cdot}$ in the equation \eqref{eq:esquema_DG_NS-CH_u}, using $\shd{\cdot}{\cdot}{\cdot}{\cdot}$ defined in \eqref{def:regularization_sh} for a certain $\delta>0$,
	the following discrete energy law holds:
	\begin{align}
	\label{discrete_energy_law_NS-CH_regularized_sh}
	\delta_tE(\uu^{m+1},\Pih_1\phi^{m+1})&+2\escalarLd{\eta(\phi^{m+1})\DD\uu^{m+1}}{\DD\uu^{m+1}}+\bupw{-\nabla_{\nn}^0\mu^{m+1}}{M(\phi^{m+1})}{\Pi_0\mu^{m+1}}\nonumber\\&\le-\frac{\delta}{2}\sum_{e\in\Ehi}\int_e\frac{\vert\uu^{m+1}\cdot\nn_e\vert}{\vert \uu^{m+1}\cdot\nn_e\vert+\delta}\salto{\Pi_0\mu^{m+1}}\salto{\phi^{m+1}}.
	\end{align}
\end{theorem}

\section{Numerical experiments}
\label{sec:numerical_experiments}

We have carried out the following numerical experiments in the spatial domain $\Omega=[-0.5,0.5]^2$. Moreover, we have set the following values of the parameters $\varepsilon=0.01$, $\lambda=0.01$, $\rho_1=1$ and $\rho_2=100$, unless otherwise specified. Also, we have chosen a constant viscosity, $\eta(\phi)=1$. Following the Remark~\ref{rmk:inf-sup_spaces}, we have chosen the pair of ``inf-sup'' stable spaces $(\Uh,\Ph)=((\Pb_2(\T_h)\cap H^1_0(\Omega))^d,\Pd_1(\T_h))$. Moreover, to comply with Hypothesis~\ref{hyp:mesh_n}, we have used a triangular mesh $\T_h$ resulting from halving a squared mesh using the diagonals.

To compute the approximations we have used the finite element library FEniCSx (see \cite{BasixJoss,AlnaesEtal2014,ScroggsEtal2022}) coupled with PyVista for the visualization of the results (see \cite{sullivan2019pyvista}). The source code for our implementation is hosted on GitHub\footnote{\url{https://github.com/danielacos/Papers-src}}. On the one hand, an iterative Newton solver has been used to approximate the nonlinear problem. In this sense, the modified stabilization term $\shd{\cdot}{\cdot}{\cdot}{\cdot}$ with $\delta=10^{-6}$ has been used in the scheme \eqref{esquema_DG_NS-CH} to avoid convergence issues. On the other hand, we have used the default iterative linear solver, GMRES (generalized minimal residual method), and preconditioner, computed using an incomplete LU factorization (ILU), of PETSc (see \cite{petsc-user-ref,DalcinPazKlerCosimo2011}) for solving the resulting linear systems.

\begin{remark}
	We must be careful when dealing with an ill-posed nonlinear problem if we want Newton's method to converge. To overcome this issue in the case of the approximation \eqref{esquema_DG_NS-CH}, we have added a penalty term $\xi\escalarL{p^{m+1}}{\overline{p}}$ to the LHS of \eqref{eq:esquema_DG_NS-CH_p} with $\xi>0$ very small (in practice, we have chosen $\xi=10^{-10}$). In this way, we enforce the $0$-mean constraint on the approximation of the pressure  $p$ and Newton's method does converge. In fact, a posteriori, we can check that this additional term has not severely affected the approximation obtained in two different manners. On the one hand, taking into account the $\norma{\cdot}_\infty$ of the approximation of $p$ we observe that the term $\xi p$ has been at most of order $10^{-5}$. On the other hand, the pointwise bounds have been preserved despite the crucial role that the local incompressibility constraint \eqref{local_incompressibility} plays in Theorem~\ref{thm:discrete_maximum_principle}.

	Certainly, many other ways of enforcing the $0$-mean pressure constraint in the nonlinear system can be explored. For instance, another interesting possibility could be adding a penalty term $\gamma\int_\Omega p^{m+1}\int_\Omega\overline{p}$, with $\gamma>0$, to the LHS of \eqref{eq:esquema_DG_NS-CH_p} as done in \cite{pacheco2023optimal}.
\end{remark}

In all the figures shown in this section, we plot both the phase field variable (in red/blue) and the following scaled vector field (in white)
$$
\uu^{m+1}_s=
\begin{cases}
	\frac{5\cdot 10^{-2}}{\norma{\uu^{m+1}}_{L^\infty(\Omega)}}\uu^{m+1},& \text{if } \norma{\uu^{m+1}}_{L^\infty(\Omega)}\ge 5\cdot 10^{-2},\\
	\uu^{m+1}, &\text{otherwise}.
\end{cases}
$$

As a reference, the computational time for these tests in a personal computer with Intel Core i7-6700 3.40GHz using 8 threads has been the following: 10 hours to compute the reference solution in Test~\ref{test:accuracy}, around 1.5 hours for Test~\ref{test:circle}, around 24 hours for Test~\ref{test:bubble} and around 33 hours for Test~\ref{test:rayleigh}.

\subsection{Accuracy test}
\label{test:accuracy}

\begin{figure}[h]
	\centering
	\includegraphics[scale=0.22]{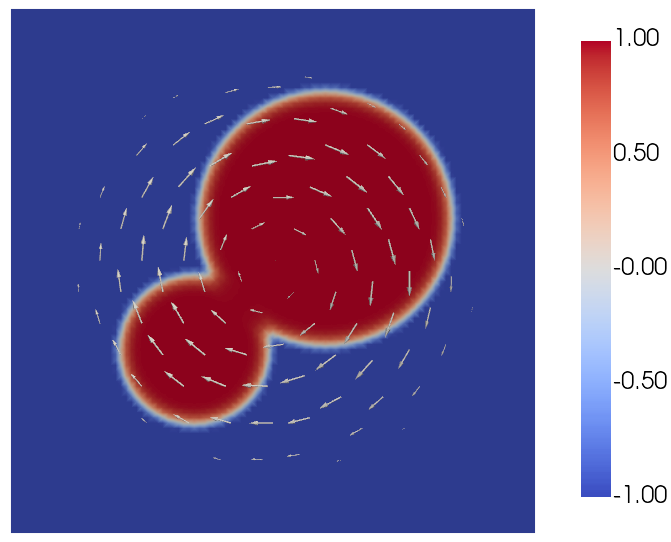}
	\caption{Initial condition of Tests~\ref{test:accuracy} and~\ref{test:circle}.}
	\label{fig:circle_init}
\end{figure}

In this case,
we define the following initial conditions
\begin{align*}
\phi_0(x,y)&=2 \tanh\left(\frac{(0.25 - \sqrt{(x-0.1)^2 + (y-0.1)^2})_\oplus}{\sqrt{2}\varepsilon}\right.\\ &\quad+ \left.\frac{(0.15 - \sqrt{(x+0.15)^2 + (y+0.15)^2})_\oplus}{\sqrt{2}\varepsilon}\right) - 1.0,\\
\uu_0(x,y)&=\chi(y(0.16-(x^2+y^2))_\oplus, -x(0.16-(x^2+y^2))_\oplus),
\end{align*}
with $\chi=1$, which are plotted in Figure~\ref{fig:circle_init}.

We conduct a preliminary convergence test in which we compare a reference solution in a very refined mesh ($h\approx 7\cdot 10^{-3}$, $962404$ degrees of freedom) with the approximation in a less refined mesh. In this way, with $\Delta t=10^{-5}$ fixed, we can remove the error introduced by the time discretization in each of the different schemes. In any case, we would like to emphasize that such a test for these sophisticated schemes involving several different discrete spaces and projection operators is nontrivial and the results obtained only provide an estimation of the possible order of convergence of the proposed approximations.

The results of the test at $T=5\cdot 10^{-4}$ are shown in Tables~\ref{table:errors_L2} and \ref{table:errors_H1}. It is worth mentioning that, as in \cite{acosta-soba_upwind_2022} for the convective Cahn-Hilliard model, order 2 in $\norma{\cdot}_{L^2(\Omega)}$ and order 1 in $\norma{\cdot}_{H^1(\Omega)}$ for the approximation of the variable $\Pih_1\phi$ have been approached. On the other hand, order around 2 in $\norma{\cdot}_{L^2(\Omega)}$ has been obtained for the approximations of $p$ and $\uu$, the latter probably affected by the order of convergence in the approximation of $\Pih_1\phi$. Finally, order around $2$ in $\norma{\cdot}_{H^1(\Omega)}$ seems to have been achieved by the approximation of $\uu$.

\begin{table}
	\centering
	\small
	\begin{tabular}{||c|c|c|c|c|c|c|c||}
		\hline
		\multirow{2}{*}{Variable} & $h\approx 2.36\cdot 10^{-2}$ & \multicolumn{2}{c|}{$3h/4\approx 1.77\cdot 10^{-2}$} & \multicolumn{2}{c|}{$4h/7\approx 1.35\cdot 10^{-2}$} &\multicolumn{2}{c||}{$h/2\approx 1.18\cdot 10^{-2}$} \\
		\cline{2-8}
		&Error &Error & \textbf{Order} &  Error & \textbf{Order} &  Error & \textbf{Order}  \\
		\hline
		\hline
		 $\Pih_1\phi$ & $8.48e-03$ & $5.40e-03$ & $1.57$ & $3.38e-03$ & $1.73$ & $2.62e-03$ & $1.89$ \\
		\hline
		\hline
		$\uu$ & $5.91e-04$ & $4.89e-04$ & $0.66$ & $3.31e-04$ & $1.44$ & $2.43e-04$ & $2.30$ \\
		\hline
		\hline
		$p$ & $2.24e-01$ & $1.14e-01$ & $2.35$ & $5.47e-02$ & $2.71$ & $4.37e-02$ & $1.67$ \\
		\hline
	\end{tabular}
	\captionof{table}{Errors and convergence orders at $T=5\cdot 10^{-4}$ in $\norma{\cdot}_{L^2(\Omega)}$.}
	\label{table:errors_L2}
\end{table}

\begin{table}
	\centering
	\small
	\begin{tabular}{||c|c|c|c|c|c|c|c||}
		\hline
		\multirow{2}{*}{Variable} & $h\approx 2.36\cdot 10^{-2}$ & \multicolumn{2}{c|}{$3h/4\approx 1.77\cdot 10^{-2}$} & \multicolumn{2}{c|}{$3h/5\approx 1.41\cdot 10^{-2}$} &\multicolumn{2}{c||}{$h/2\approx 1.18\cdot 10^{-2}$} \\
		\cline{2-8}
		&Error &Error & \textbf{Order} &  Error & \textbf{Order} &  Error & \textbf{Order}  \\
		\hline
		\hline
        $\Pih_1\phi$ &  $1.22e+00$ & $1.17e+00$ & $0.15$ & $9.12e-01$ & $0.92$ & $8.09e-01$ & $0.89$ \\
		\hline
		\hline
		 $\uu$ & $9.61e-02$ & $7.98e-02$ & $0.65$ & $4.90e-02$ & $1.80$ & $3.75e-02$ & $1.99$ \\
		\hline
	\end{tabular}
	\captionof{table}{Errors and convergence orders at $T=5\cdot 10^{-4}$ in $\norma{\cdot}_{H^1(\Omega)}$.}
	\label{table:errors_H1}
\end{table}

\begin{remark}
	Several works such as \cite{diegel2017convergence,chen2022error,chen2022errorCHNS,styles2008finite} have carried out a careful error analysis of finite element approximations of phase-field models coupled with fluid motion such as the CHNS system or related models. However, most of these works have focused on the case of constant or non-degenerate mobility and constant density and their results are based on the energy-stability property of the proposed approximations. It is left for a future work to study whether these techniques can be extended and applied to derive error estimates for our proposed approximation \eqref{esquema_DG_NS-CH}.
\end{remark}

\subsection{Mixing bubbles}
\label{test:circle}

For this test we keep the same initial conditions as in the previous test but with $\chi=100$. Again, this initial condition can be seen in Figure~\ref{fig:circle_init}.

\begin{figure}[H]
	\centering
	\begin{tabular}{ccc}
		 \hspace*{-1cm}$t=2\cdot 10^{-2}$ & \hspace*{-1cm}$t=5\cdot 10^{-2}$ & \hspace*{-1cm}$t=10^{-1}$ \\
		\includegraphics[scale=0.22]{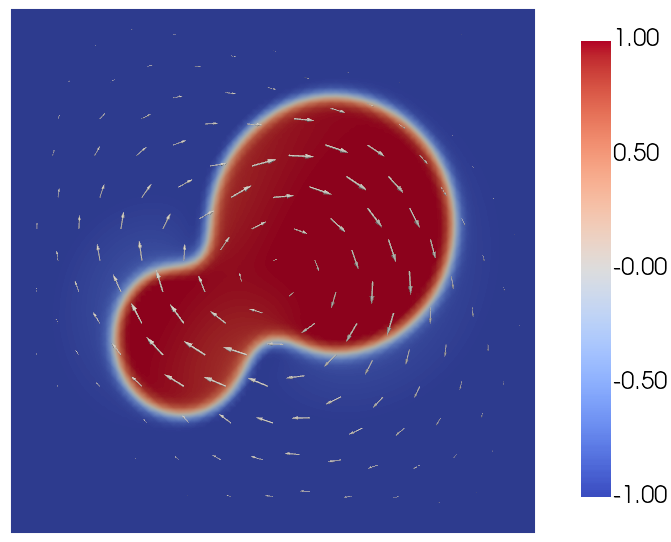} &
		\includegraphics[scale=0.22]{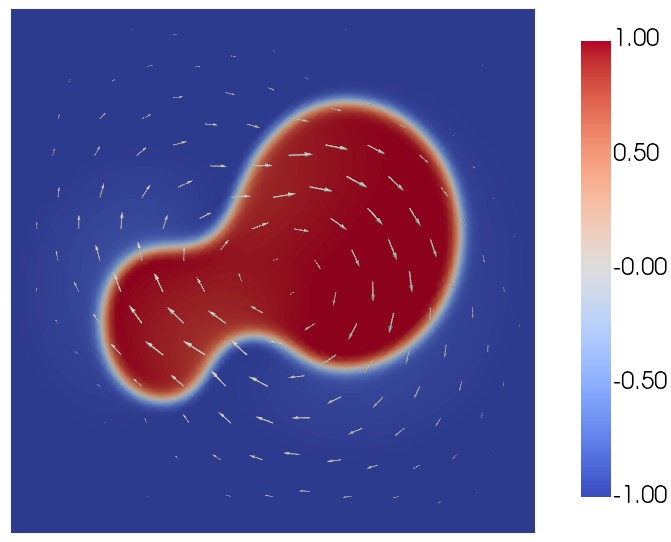} &
		\includegraphics[scale=0.22]{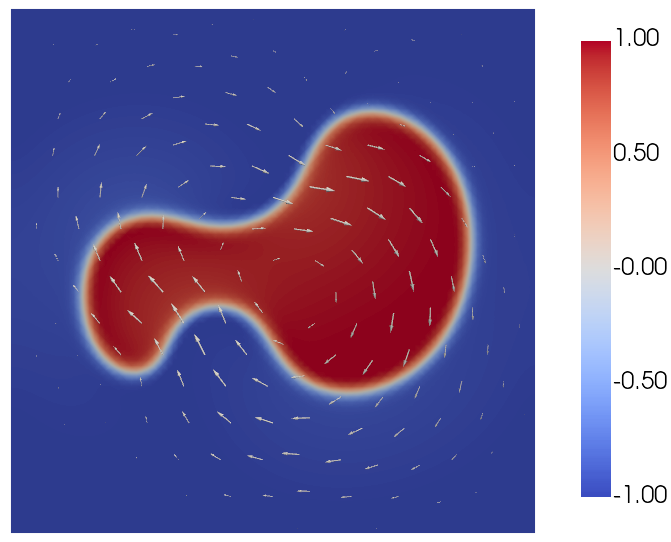} \\
	\end{tabular}
	\caption{Evolution of $\Pih_1\phi$ over time in Test~\ref{test:circle}.}
	\label{fig:circle_evol}
\end{figure}

\begin{figure}[H]
	\centering
	\begin{tabular}{cc}
		\includegraphics[scale=0.52]{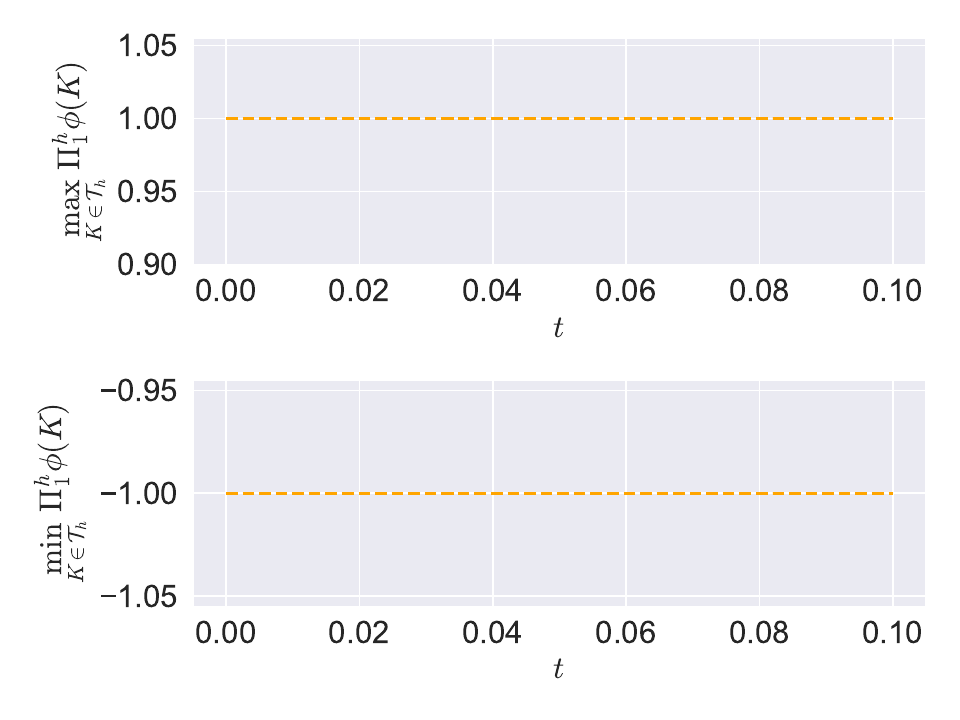} &
		\includegraphics[scale=0.52]{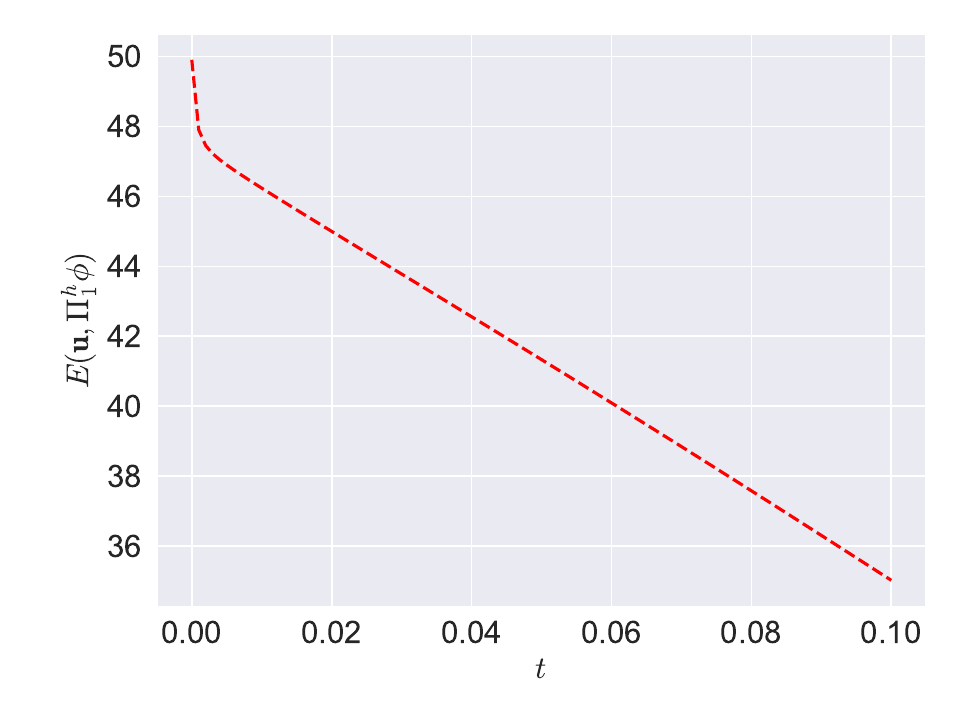}
	\end{tabular}
	\caption{Left, maximum and minimum of $\Pih_1\phi$. Right, discrete energy. Test~\ref{test:circle}.}
	\label{fig:circle_max-min_energy}
\end{figure}

In Figure~\ref{fig:circle_evol} we have plotted the evolution in time of the approximation obtained using both the scheme \eqref{esquema_DG_NS-CH} with $h\approx 1.41\cdot 10^{-2}$ ($241204$ degrees of freedom) and $\Delta t=10^{-3}$. On the other hand, in Figure~\ref{fig:circle_max-min_energy} (left) we can observe how the bounds are preserved as predicted by the previous analytical results. In addition, in Figure~\ref{fig:circle_max-min_energy} (right) one may observe how the energy decreases as predicted by the theory above.

\begin{figure}[H]
	\centering
	\begin{tabular}{ccc}
		\hspace*{-1cm}$t=0$ & \hspace*{-1cm}$t=6.5\cdot10^{-3}$ & \hspace*{-1cm}$t=1.2\cdot 10^{-2}$\\
		\includegraphics[scale=0.22]{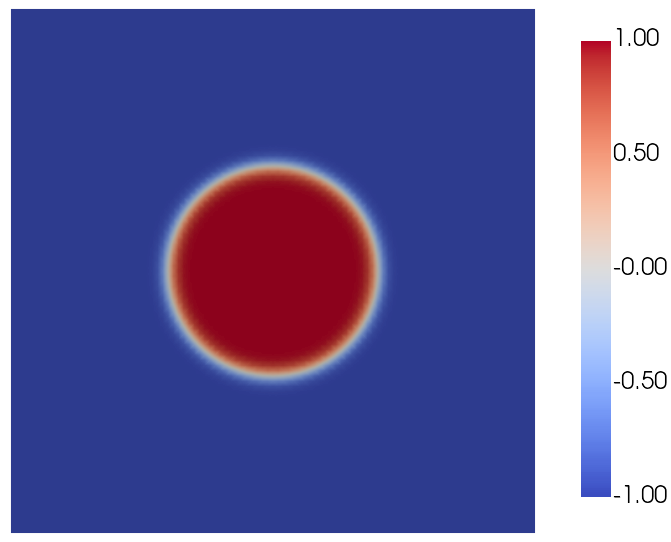} &
		\includegraphics[scale=0.22]{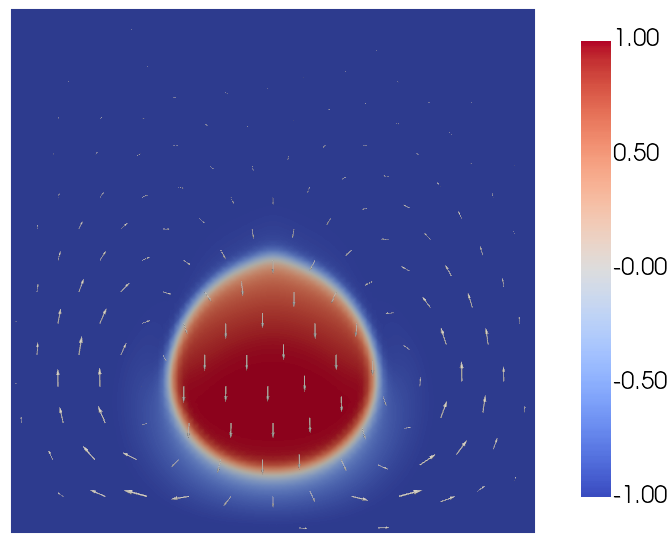} &
		\includegraphics[scale=0.22]{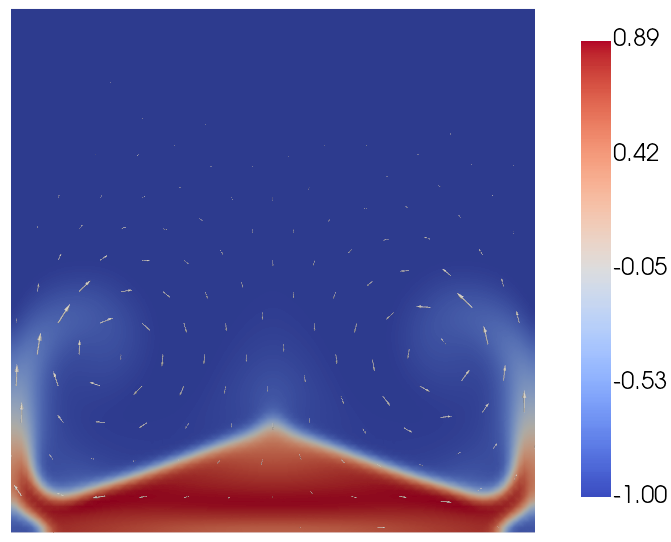} \\
		\hspace*{-1cm}$t=3\cdot 10^{-2}$ & \hspace*{-1cm}$t=4.5\cdot 10^{-2}$ & \hspace*{-1cm}$t=2.5\cdot 10^{-1}$ \\ \includegraphics[scale=0.22]{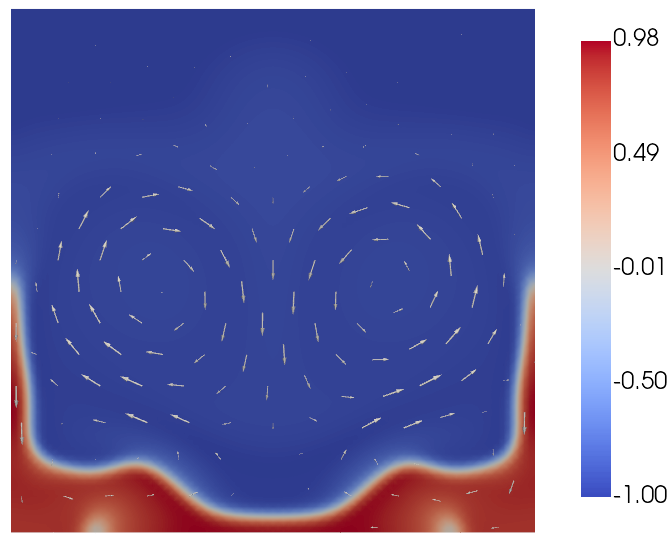} &
		\includegraphics[scale=0.22]{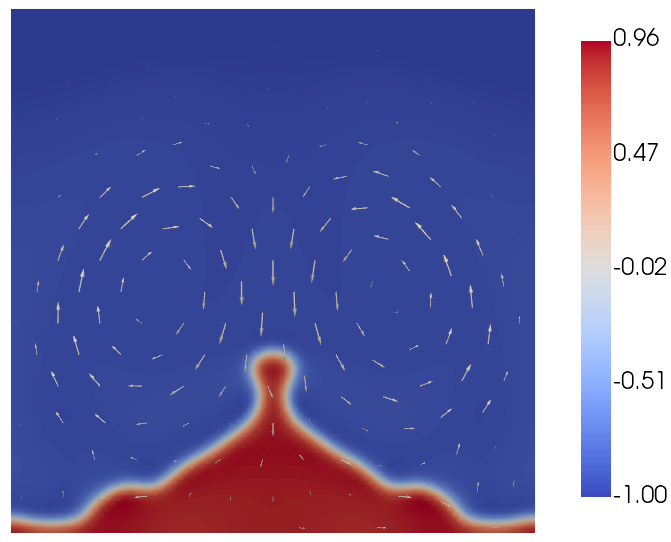} &
		\includegraphics[scale=0.22]{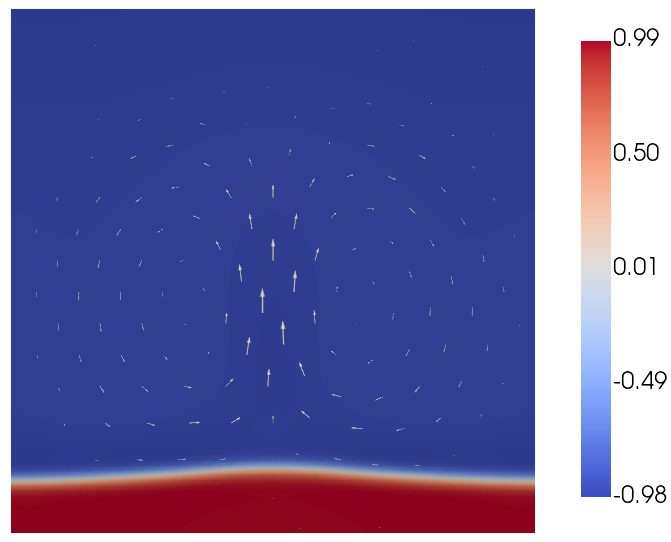}
	\end{tabular}
	\caption{Evolution of $\Pih_1\phi$ over time in Test~\ref{test:bubble}.}
	\label{fig:bubble_evol}
\end{figure}

\subsection{A heavier bubble falling in a lighter medium}
\label{test:bubble}

Now, we perform a test in which we define the following initial condition: $\uu_0=0$ and 
$$
\phi_0(x,y)=\tanh\left(\frac{0.2 - \sqrt{x^2 + y^2}}{\sqrt{2}\varepsilon}\right),
$$
a bubble of density $\rho_2=100$ in a lighter medium of density $\rho_1=1$, plotted in Figure~\ref{fig:bubble_evol} ($t=0$). Moreover, we have added a term $-\rho(\phi)\boldsymbol{g}$ on the right-hand side of equation \eqref{eq:NS-CH_model_u} acting as the gravitational forces pushing the heavier bubble down to the bottom of the domain $\Omega$. In our case, we have chosen $\boldsymbol{g}=(0,1)$ and we have treated this term implicitly in \eqref{esquema_DG_NS-CH}.

\begin{figure}[H]
	\centering
	\begin{tabular}{cc}
		\includegraphics[scale=0.52]{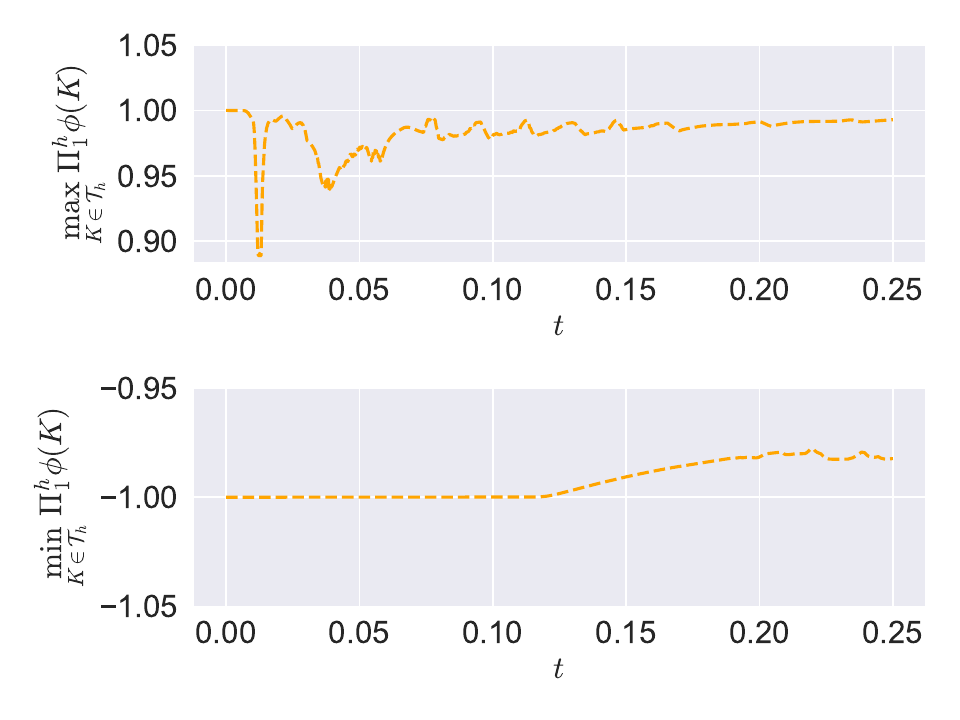} &
		\includegraphics[scale=0.52]{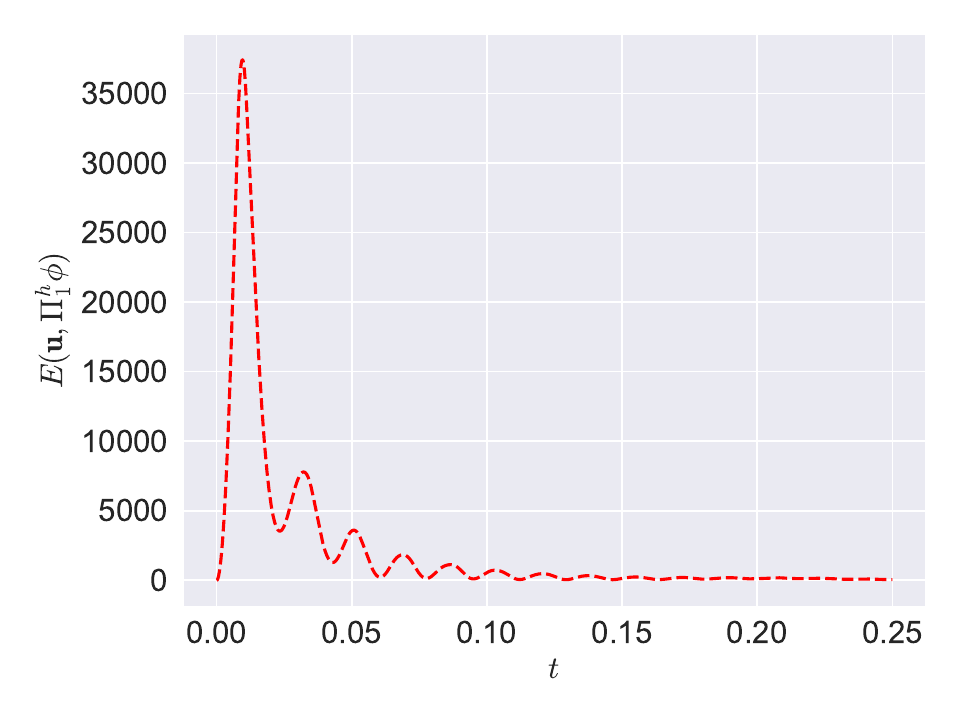}
	\end{tabular}
	\caption{Left, maximum and minimum of $\Pih_1\phi$. Right, discrete energy. Test~\ref{test:bubble}.}
	\label{fig:bubble_max-min_energy}
\end{figure}

In this case, we have shown in Figure~\ref{fig:bubble_evol} the evolution in time of the solution using \eqref{esquema_DG_NS-CH} with $h\approx 1.41\cdot 10^{-2}$ and $\Delta t=10^{-4}$. The result is qualitatively similar to the ones shown in previous studies such as \cite{tierra_guillen_abels_2014}. Also, the bounds are preserved as shown in Figure~\ref{fig:bubble_max-min_energy} (left). In this case, the energy does not necessarily decrease due to the gravitational forces as one may observe in Figure~\ref{fig:bubble_max-min_energy} (right).

\subsection{A Rayleigh-Taylor type instability}
\label{test:rayleigh}

\begin{figure}
	\centering
	\begin{tabular}{ccc}
		\hspace*{-1cm}$t=0$ & \hspace*{-1cm}$t=1.25\cdot10^{-2}$ & \hspace*{-1cm}$t=2\cdot10^{-2}$\\
		\includegraphics[scale=0.22]{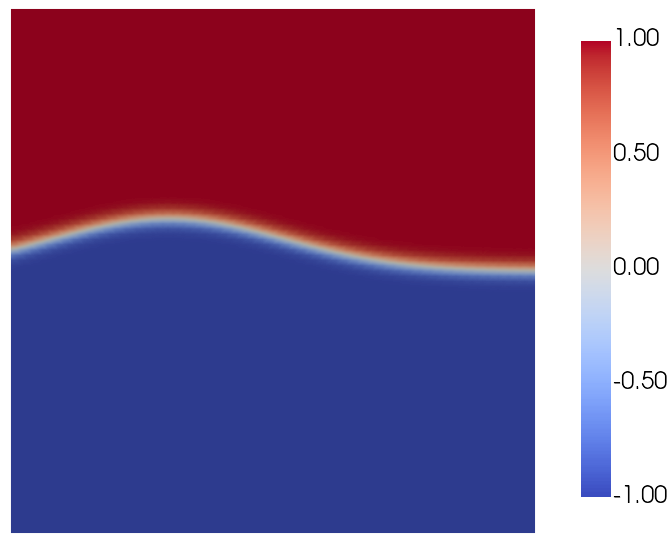} &
		\includegraphics[scale=0.22]{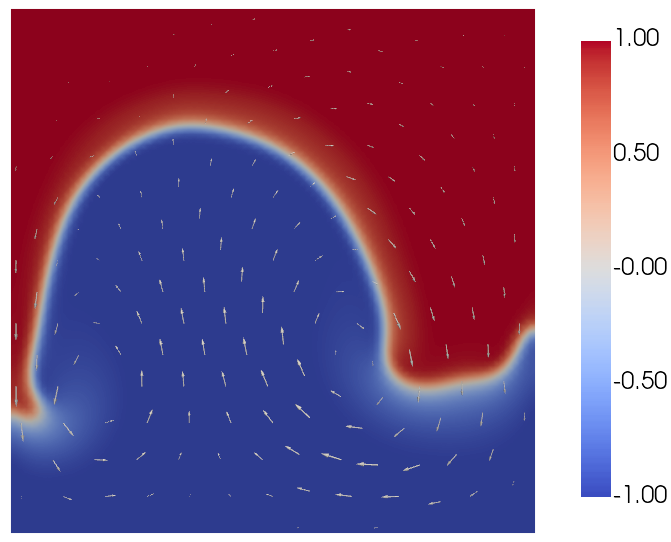} &
		\includegraphics[scale=0.22]{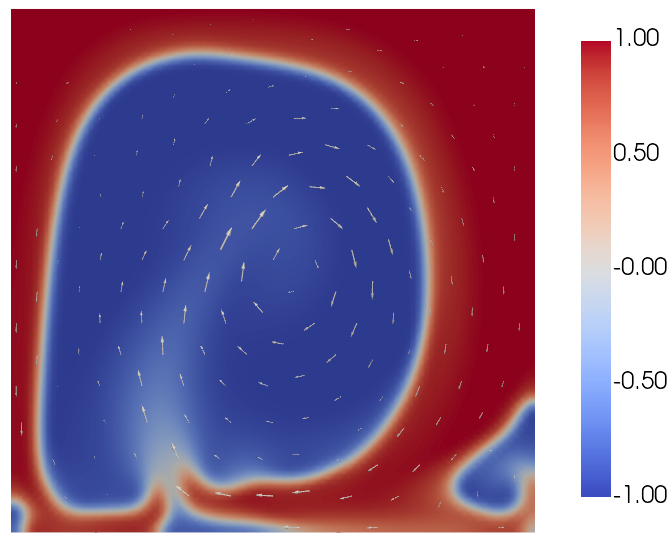} \\
		\hspace*{-1cm}$t=3\cdot10^{-2}$ & \hspace*{-1cm}$t=8\cdot10^{-2}$ & \hspace*{-1cm}$t=3.5\cdot10^{-1}$ \\
		\includegraphics[scale=0.22]{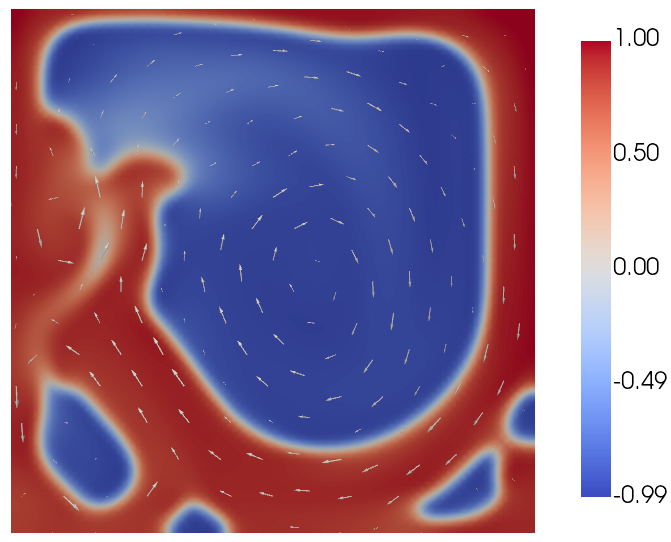} & 
		\includegraphics[scale=0.22]{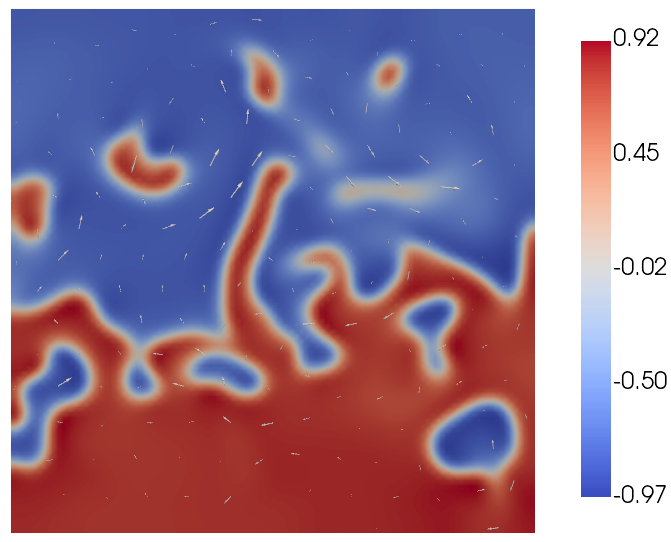} & 
		\includegraphics[scale=0.22]{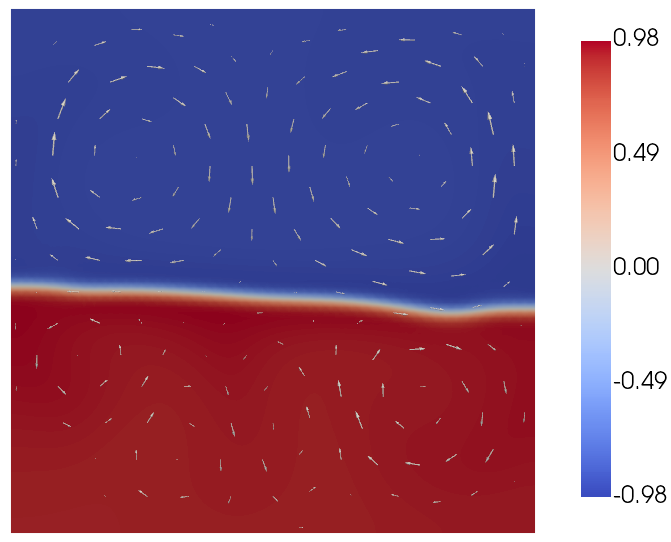}
	\end{tabular}
	\caption{Evolution of $\Pih_1\phi$ over time in Test~\ref{test:rayleigh}.}
	\label{fig:rayleigh_evol}
\end{figure}

Finally, we carry out a benchmark Rayleigh-Taylor type instability test based on the one shown in \cite{tierra_guillen_abels_2014} for which we define the following initial condition: $\uu_0=0$ and 
$$
\phi_0(x,y)=\tanh\left(\frac{y - (0.1\exp(-(x+0.2)^2/0.1))}{\sqrt{2}\varepsilon}\right),
$$
plotted in Figure~\ref{fig:rayleigh_evol} ($t=0$). Again, we add the gravity term $-\rho(\phi)\boldsymbol{g}$ with $\boldsymbol{g}=(0,1)$ in the RHS of equation \eqref{eq:NS-CH_model_u}.

\begin{figure}[H]
	\centering
	\begin{tabular}{cc}
		\includegraphics[scale=0.52]{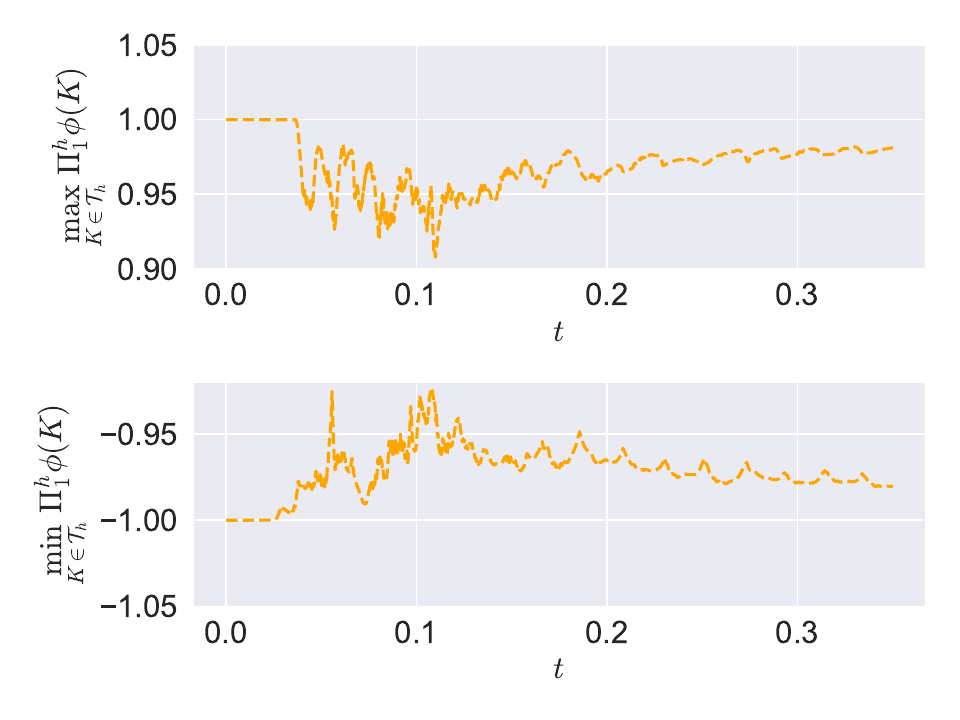} &
		\includegraphics[scale=0.52]{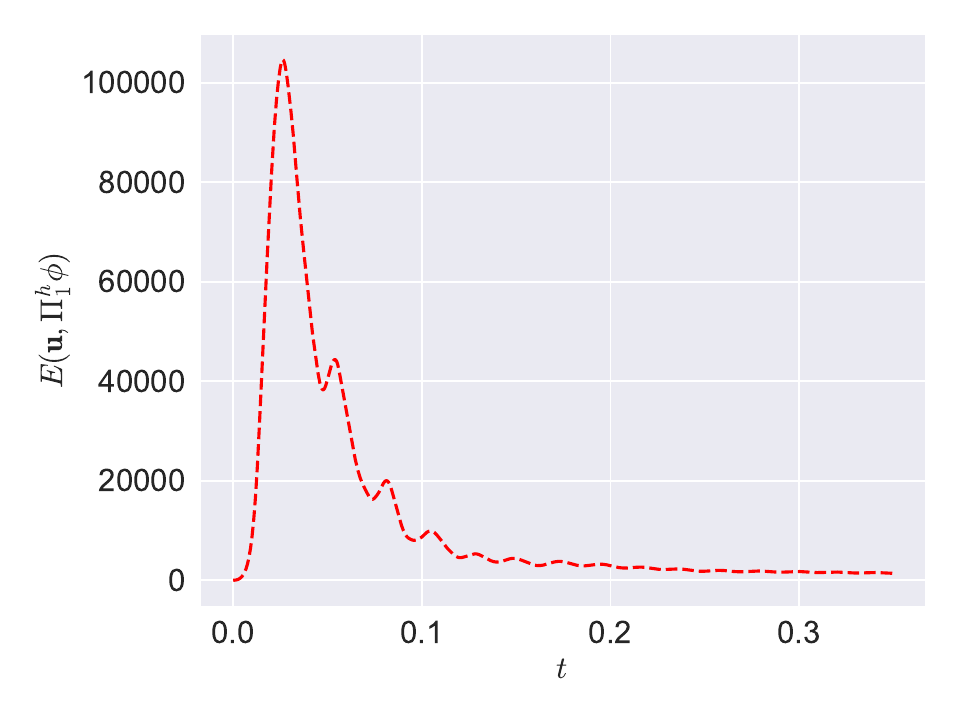}
	\end{tabular}
	\caption{Left, maximum and minimum of $\Pih_1\phi$. Right, discrete energy. Test~\ref{test:rayleigh}.}
	\label{fig:rayleigh_max-min_energy}
\end{figure}

The evolution in time of the solution using \eqref{esquema_DG_NS-CH} with $h\approx 1.41\cdot 10^{-2}$ and $\Delta t=10^{-4}$ can be seen in Figure~\ref{fig:rayleigh_evol}. Again, despite the difficulty of this test due to the fast dynamics involved, the results are qualitatively similar to the ones shown in previous works such as \cite{tierra_guillen_abels_2014}. In Figure~\ref{fig:rayleigh_max-min_energy} (left) we plot the evolution of the maximum and minimum of the regularized phase-field function, where we can observe that the bounds are indeed preserved as predicted by the theory. 
In addition, one may observe in Figure~\ref{fig:rayleigh_max-min_energy} (right) the behavior of the discrete energy.

\section{Conclusion}

In this work we have developed a robust, structure-preserving approximation, given in \eqref{esquema_DG_NS-CH}, of the CHNS model with variable density \eqref{NS-CH_model}. To our best knowledge this is the first approximation of a CHNS model with a Ginzburg-Landau polynomial potential and degenerate mobility that ensures the mass-conservation, pointwise bounds and energy-stability properties at the same time.

This approximation combines the ideas of the previous works \cite{acosta-soba_upwind_2022,acosta2023structure} to preserve the pointwise bounds of the phase-field variable as shown in Theorem~\ref{thm:discrete_maximum_principle}. In this regard, we have used a finite element approximation for the Navier-Stokes fluid flow with discontinuous pressure that preserves the incompressibility of the velocity locally in  each of the elements of the mesh $\T_h$, see \eqref{local_incompressibility}. In addition, a carefully developed upwind discontinuous Galerkin approximation for the Cahn-Hilliard part has been chosen.

Moreover, the ideas in \cite{acosta-soba_KS_2022,acosta2023structure} about approximating the normal derivative of the chemical potential in a structured mesh, \eqref{eq:approx_gradn}, and the bilinear form \eqref{def:bupw_saltos} have been employed. These ideas have been combined with novel stabilization techniques such as \eqref{def:sh2} and \eqref{upwind_stabilization_energy}, and the stabilization term \eqref{def:sh1} that was previously developed in \cite{tierra_guillen_abels_2014}. This approach has led us to the discrete energy-stability property shown in Theorem~\ref{thm:discrete_energy_law_NS-CH}.

Finally, the theoretical discussion has been complemented with several numerical experiments where the good properties of the approximation proposed are manifested. In Test~\ref{test:accuracy}, a preliminary accuracy test was carried out where second order of convergence seems to be achieved in $\norma{\cdot}_{L^2(\Omega)}$. Then, a qualitative Test~\ref{test:circle} was computed, where the discrete energy-stability property has been exhibited. Finally, two benchmark problems where the action of gravitational forces has been taken into account were conducted: a heavier bubble in a lighter medium (Test~\ref{test:bubble}) and a Rayleigh-Taylor type instability (Test~\ref{test:rayleigh}). Throughout these tests it could be seen how the pointwise bounds are preserved at the discrete level for the phase-field variable.

Despite the robustness and good properties of this new numerical approximation, we would like to mention that there is still much room for improvement. In particular, the main drawback of this numerical scheme is the computational cost inherent to such a fully coupled approximation.

In this sense, we have also explored the idea of developing a decoupled property-preserving approximation of \eqref{NS-CH_model} by means of a rotational pressure projection technique following the previous work in \cite{liu2022pressure}. However, this has been finally left for a future work due to the number of difficulties related to such kind of approximations. On the one hand, applying a rotational projection technique to a model with variable viscosity is not trivial as shown in \cite{deteix2018improving,deteix2019shear,plasman2020projection}. On the other hand, developing a stable decoupled approximation for a system involving variable densities requires carefully adjusting the intermediate steps as in \cite{pyo2007gauge,guermond2000projection,guermond2009splitting}. In addition, preserving both the pointwise bounds for the phase-field variable and the energy law of the system at the discrete level requires imposing additional restrictions, such as \eqref{local_incompressibility} and \eqref{eq:approx_gradn}, on the techniques implemented. A preliminary work on a decoupled approximation for this system \eqref{NS-CH_model}, in the case of constant viscosity, that preserves the pointwise bounds can be seen in \cite[Section 6.5]{acosta2024analysis}.

\section*{CRediT authorship contribution statement}
\textbf{Daniel Acosta-Soba:} Writing -- review \& editing, Writing -- original draft, Conceptualization, Formal Analysis, Methodology, Software, Visualization, Investigation. \textbf{Francisco Guillén-González:} Writing -- review \& editing, Supervision, Conceptualization, Formal Analysis, Methodology. \textbf{J. Rafael Rodríguez Galván:} Writing -- review \& editing, Supervision, Conceptualization, Formal Analysis, Methodology, Software. \textbf{Jin Wang:} Writing -- review \& editing, Supervision, Conceptualization, Formal Analysis, Methodology.

\section*{Acknowledgements}
The first author has been supported by \textit{UCA FPU contract UCA/REC14VPCT/2020 funded by Universidad de Cádiz} and by a \textit{Graduate Scholarship funded by The University of Tennessee at Chattanooga}. The second and third authors have been supported by 
 \textit{Grant PGC2018-098308-B-I00 (MCI/AEI/FEDER, UE, Spain)}, \textit{Grant US-1381261 (US/JUNTA/FEDER, UE, Spain)}
and \textit{Grant P20-01120 (PAIDI/JUNTA/FEDER, UE, Spain)}.
The fourth author has been supported by the US National Science Foundation 
under Grant Numbers 1913180 and 2324691. 

Furthermore, we would like to thank the reviewers for their thoughtful and mind-opening comments that have helped us significantly improve our work. In particular, we are grateful for their guidance on fractional time-stepping approximations that have shed light on possible future works.

\printbibliography

\end{document}